\documentclass[11pt]{article}
\usepackage[a4paper, margin = 2.2cm]{geometry}
\usepackage[english]{babel}
\usepackage{setspace}
\usepackage{latexsym,amsmath,enumerate,graphics,enumerate,amsthm,tikz,hyperref,float}
\usepackage[mathscr]{euscript}
\usepackage[affil-it]{authblk}
\usepackage{enumerate,amsthm,dsfont,pstricks}
\usepackage{latexsym,amsmath,amssymb,amscd,wrapfig,graphicx}
\usepackage{wrapfig,graphicx,curves}
\usepackage{sectsty}
\sectionfont{\centering}

\newtheorem{theorem}{Theorem}[section]
\newtheorem{lemma}[theorem]{Lemma}
\newtheorem{proposition}[theorem]{Proposition}
\newtheorem{corollary}[theorem]{Corollary}

\theoremstyle{definition}
\newtheorem{definition}[theorem]{Definition}

\newtheorem{remark}[theorem]{Remark}
\theoremstyle{remark}

\let\phi=\varphi
\def\epsilon{\varepsilon}

\def\0{\mathbf{0}}

\DeclareMathOperator{\diam}{diam}

\newcommand{\comment}[1]{}

\numberwithin{equation}{section}

\let\epsilon=\varepsilon

\makeatletter
\def\@maketitle{%
  \newpage
  \null
  \vskip 2em%
  \begin{center}%
  \let \footnote \thanks
    {\Large\bfseries \@title \par}%
    \vskip 2.5em%
    {\normalsize
      \lineskip .5em%
      \begin{tabular}[t]{c}%
        \@author
      \end{tabular}\par}%
   \vskip 1em%
  \end{center}%
  \par
  \vskip 1em}
\makeatother

\begin{document}

\title{\bf  \huge Horofunction compactifications of symmetric cones under Finsler distances}

\author{Bas Lemmens%
\thanks{Email: \texttt{B.Lemmens@kent.ac.uk}, Bas Lemmens gratefully acknowledges the support of the EPSRC (grant EP/R044228/1)}}
\affil{School of Mathematics, Statistics \& Actuarial Science,
University of Kent, \\Canterbury, CT2 7NX, United Kingdom\\
\mbox{}\\
}

\maketitle
\date{}

\begin{abstract}
In this paper we consider symmetric cones as symmetric spaces equipped with invariant Finsler distances, namely the Thompson distance and the Hilbert distance. We establish a correspondence between the horofunction compactification of a symmetric cone $A_+^\circ$  under these invariant Finsler distances and the horofunction compactification of the normed space in the tangent bundle.  More precisely, for  the Thompson distance on $A^\circ_+$ we show that the exponential map extends as a homeomorphism between the horofunction compactification of the normed space in the tangent bundle, which is a JB-algebra, and the horofunction compactification of $A_+^\circ$.  We give a complete  characteristation of the Thompson distance horofunctions and provide an explicit extension of the exponential map.  Analogues results are established for the Hilbert distance on the projective cone $PA_+^\circ$. The analysis yields a geometric description of the horofunction compactifications of these spaces in terms of the facial structure of the closed unit ball of the dual norm of the norm in the tangent space. 
\end{abstract}



\section{Introduction}
Symmetric cones provide an important class of Riemannian symmetric spaces of non-compact type. They support many invariant Finsler metrics, that are not Riemannian, which turn them into  Finsler symmetric spaces. In particular, every symmetric cone can be identified as the interior, $A_+^\circ$, of the cone of  squares  in a finite dimensional formally real Jordan algebra $A$ by the Koecher-Vinberg theorem.   The  real Jordan algebra $A$ with unit $u$ can be equipped with the spectral norm: 
\[
\|x\|_u =\inf\{\mu>0\colon -\mu u\leq x\leq \mu u\} =\max\{|\lambda|\colon \lambda\in \sigma(x)\},
\] 
where $x\leq y$ if $y-x$ is positive, and $\sigma(x)$ is the spectrum of $x$, which turns it into a JB-algebra \cite{AS1}. One can use the spectral norm to  put an invariant Finsler distance on $A_+^\circ$ by defining the Finsler metric $F\colon TA_+^\circ\to \mathbb{R}$  by 
$F(x,w) = \|U_{x^{-1/2}}w\|_u$, where $U_z\colon A\to A$ is the quadratic representation of $z$. (It should be noted that in this setting $F$ does not satisfy the smoothness and strong convexity conditions commonly used in the theory of Riemann-Finsler manifolds \cite{BCS}.) The Finsler distance between $x$ and $y$ in $A_+^\circ$ is the infimum of lengths, 
\[
L(\gamma) = \int_0^1 F(\gamma(t),\gamma'(t)) \,\mathrm{d} t,
\] 
 over all piecewise $C^1$-smooth paths $\gamma\colon [0,1]\to A_+^\circ$ with $\gamma(0)=x$ and $\gamma(1)=1$.  
In fact, this distance is known \cite{Lim,Nu} to coincide with the Thompson distance on $A^\circ_+$, which is given by  
\[
d_T(x,y) = \max\{|\mu|\colon \mu\in \sigma(\log U_{y^{-1/2}}x)\}  =   \|\log U_{x^{-1/2}}y\|_u  
\]
for all $x,y\in A_+^\circ$. 

In this setting the Riemannian symmetries $S_x\colon A^\circ_+\to A_+^\circ$, which are given by $S_x(y) = U_{x} y^{-1}$, are global $d_T$-isometries (see \cite{LRW})  that make $A_+^\circ$  a  Finsler symmetric space. Moreover, the automorphism group $\mathrm{Aut}(A_+)=\{T\in\mathrm{GL}(A)\colon T(A_+)=A_+\}$ is a subgroup of the isometry group of $(A_+^\circ,d_T)$ and acts transitively on $A_+^\circ$. 

A prime example is the symmetric cone of $n\times n$ strictly positive definite Hermitian matrices $\Pi_n(\mathbb{C})$, which can  be  identified with the symmetric space $\mathrm{GL}_n(\mathbb{C})/\mathrm{U}_n$ by letting $\mathrm{GL}_n(\mathbb{C})$ act  on $\Pi_n(\mathbb{C})$ by $A\mapsto MAM^*$.  In this setting the quadratic representation $U_A\colon B\mapsto ABA$, so the Finsler metric in the tangent space at $A$ satisfies  $F(A,M) = \|A^{-1/2}MA^{-1/2}\|_I = \max\{|\lambda|\colon \lambda \in \sigma(A^{-1/2}MA^{-1/2})\}$, and   the Thompson distance is given by $d_T(A,B) = \|\log A^{-1/2}BA^{-1/2}\|_I = \max\{|\lambda|\colon \lambda \in \sigma(\log A^{-1/2}BA^{-1/2})\}$. 

Given a symmetric cone $A_+^\circ$ it also interesting to consider the projective cone $PA_+^\circ = A_+^\circ/\mathbb{R}_{>0}$, so $x\sim \lambda x$ for $x\in A_+^\circ$ and $\lambda>0$. We will identify this space with the set of points $x$ in $A_+^\circ$ with $\det x =1$, that is, 
\[
PA_+^\circ=\{x\in A_+^\circ\colon \det x =1\}.
\] 
The tangent space at the unit $u\in PA_+^\circ$ is given by $T_u = \{w\in A\colon \mathrm{tr}\, w =0\}$, as $(D_w\det)(x)=(\det x)\mathrm{tr}(U_{x^{-1/2}}w)$ for $x\in A_+^\circ$, see \cite[p.\,53]{FK}. In fact, the tangent space at $x\in PA_+^\circ$ is $T_x= U_{x^{1/2}} (T_u)$. 

There is a natural Finsler metric $H$ on the tangent bundle. For the unit $u\in PA_+^\circ$ and $w\in T_u$,  let 
\[
H(u,w)  = \inf\{\lambda \colon w\leq \lambda u\} -\sup\{\lambda \colon  \lambda u\leq w\} = \max\sigma(w)-\min\sigma(w)  =\mathrm{diam}\,\sigma(w).
\]
Note that $H(u,\cdot)$ is a norm on $T_u$, which  is called the variation norm and will be denoted by $|\cdot|_u$. For $x\in PA_+^\circ$ and $w\in T_x$, let 
$H(x,w) = |U_{x^{-1/2}} w|_u$. 
 
In this case the Finsler distance $d_H$ on $PA_+^\circ$ coincides with the Hilbert distance on $PA^\circ_+$, see \cite{Nu} or \cite[Proposition 5.3]{LP},  and satisfies 
\[
d_H(x,y) = \max \sigma(\log U_{x^{-1/2}}y)  - \min \sigma(\log U_{x^{-1/2}}y) =  \mathrm{diam}\, \sigma(\log U_{x^{-1/2}}y). 
\]
for all $x,y\in PA_+^\circ$.  The space $PA_+^\circ$ is a Finsler symmetric space, where the  Riemannian symmetry $S_x\colon y \mapsto U_x y^{-1}$ for $y\in PA_+^\circ$, is a global $d_H$-isometry, see \cite{LRW}. Note that $S_x$ maps $PA_+^\circ$ into itself, as $\det U_x y^{-1} = (\det x)^2(\det y)^{-1}$, see \cite[Proposition III.4.2]{FK}.  For the space $P\Pi_n(\mathbb{C})$, which corresponds to the symmetric space $\mathrm{SL}_n(\mathbb{C})/\mathrm{SU}_n$, we have that $d_H(A,B) = \mathrm{diam}\,\sigma(\log A^{-1/2}BA^{-1/2})$. 
  
In this paper we consider the horofunction compactifications of the Finsler symmetric spaces  $(A_+^\circ,d_T)$  and  $(PA_+^\circ,d_H)$. We show that the exponential map extends as a homeomorphism  between the  horofunction compactifications of  the JB-algebra $(A,\|\cdot\|_u)$, i.e., the tangent space at $u$, and $(A_+^\circ, d_T)$. In addition, we establish a similar result for  the horofunction compactifications of  $(PA_+^\circ,d_H)$ and  the normed space $(T_u, |\cdot|_u)$.  More precisely, we prove the following  theorem. 
\begin{theorem}\label{main} Let $A_+^\circ$ be a symmetric cone in a finite dimensional formally real Jordan algebra $A$. 
\begin{enumerate}[(a)]
\item The exponential map, $ \mathrm{exp}_u\colon x\in A\mapsto e^x\in A_+^\circ$, extends as a homeomorphism between the horofunction compactification of the JB-algebra $(A,\|\cdot\|_u)$  and the horofunction compactification of $(A_+^\circ,d_T)$, such that each  part  in the horofunction boundary of $A$  is mapped  onto a part of the horofunction boundary of $A_+^\circ$. 
\item The exponential map, $ \mathrm{exp}_u\colon x\in T_u\mapsto e^x\in PA_+^\circ$, extends as a homeomorphism between the horofunction compactification of  $(T_u,|\cdot|_u)$  and the horofunction compactification of $(PA_+^\circ,d_H)$, such that each  part  in the horofunction boundary of $T_u$  is mapped  onto a part of the horofunction boundary of $PA_+^\circ$. 
\end{enumerate}
\end{theorem}  
To prove  Theorem \ref{main} we give complete descriptions of  the horofunction compactifications of $(A_+^\circ,d_T)$ and $(T_u,|\cdot|_u)$, together with their parts and the detour distances, see Theorems \ref{horofunctions} and \ref{detour}, Theorem \ref{horoTu} and Propositions \ref{detourTu} and \ref{partsTu}. Combining these theorems with the results in \cite{LP} allows us to give explicit extensions of the exponential map to the horofunction boundaries.

The horofunction compactification of a finite dimensional JB-algebra $(A,\|\cdot\|_u)$ was determined in \cite[Section 4]{LP}, and it was shown there that it is naturally homeomorphic to the closed unit ball  in the dual space of $(A,\|\cdot\|_u)$.  Combining the results from \cite[Section 4]{LP} with Theorem \ref{main}(a) gives the following geometric description of the horofunction compactification of $(A_+^\circ,d_T)$. 
\begin{corollary} There exists a homeomorphism between the horofunction compactification of $(A_+^\circ,d_T)$ onto the closed unit ball in the dual space of $(A,\|\cdot\|_u)$, which maps each part of the horofunction boundary onto a relative open boundary face of the ball. 
\end{corollary}

With regard to part (b) of Theorem \ref{main} it is known that there exists a homeomorphism between  the horofunction compactification of $(PA_+^\circ,d_H)$ and the  closed unit ball in the dual space of $(T_u,|\cdot|_u)$, see \cite[Section 5]{LP}. So Theorem \ref{main}(b) has  the following corollary. 
\begin{corollary}\label{cor2} There exists a homeomorphism between the horofunction compactification of $(T_u,|\cdot|_u)$ onto the closed unit ball  in the dual space of $(T_u,|\cdot |_u)$, which maps each part of the horofunction boundary onto a relative open boundary face of the ball. 
\end{corollary}

The results in this paper are motivated  by work of Kapovich and Leeb \cite{KL} and subsequent works by, Ji and Schilling \cite{JS1,JS}, Haettel, Schilling, Walsh and Wienhard  \cite{HSWW}, Schilling \cite{Sc}, and the author and Power \cite{LP}. In \cite[Question 6.18]{KL} the question was raised if the horofunction compactification of a finite dimensional normed space is naturally  homeomorphic to the closed unit ball of the dual normed space. Corollary \ref{cor2} shows that this is the case for the normed spaces $(T_u,|\cdot|_u)$. More generally, it is interesting to understand  when there exists a natural homeomorphism between the horofunction compactification of a Finsler symmetric space and the closed unit ball of the dual norm of the norm in the tangent space. It is also interesting to determine when  there exists a homeomorphism between the horofunction compactification of a  Finsler symmetric space and the horofunction compactification of the normed space in the tangent bundle at the basepoint.  Theorem \ref{main} provides  interesting classes of Finsler symmetric spaces for which both questions have a positive answer.  

Compactifications of symmetric spaces is a rich subject, which has been studied extensively, see  for instance \cite{BJ,GJT}. In recent years   Finsler structures  have been used to study Satake compactifications of symmetric spaces. In particular, Haettel, Schilling, Walsh and Wienhard  \cite{HSWW}, see also \cite{Sc}, showed that each generalised Satake compactification can be realised as a horofunction compactification with respect to an invariant  polyhedral Finsler metric on the flats. Kapovich and Leeb  \cite{KL} realised the maximal Satake compactification of a symmetric space $X=G/K$ of non-compact type as the horofunction compactification with respect to an invariant Finsler metric on $X$.  
Friedland and Freitas \cite{FF2} showed that the horofunction compactification of the Siegel upper half plane of rank $n$ under the  Finsler $1$-metric agrees with the bounded symmetric domain compactification, which is a minimal Satake compactification.  Greenfield and Ji \cite{GJ} used invariant Finsler distances  to study compactifications of the Teichm\"uller spaces of flat tori. In particular, they study the horofunction compactification of $\mathrm{SL}_n(\mathbb{R})/\mathrm{SO}_n$ under the Thompson distance, which they refer to as the generalised Hilbert metric. For the Satake and Martin compactifications of symmetric spaces  Ji \cite{Ji} showed that they are homeomorphic to the closed unit ball in the tangent space, see also \cite{Ku}.  

The Hilbert and Thompson distances are important invariant Finsler distances on the classes of symmetric spaces we consider here. It would be interesting to know if their horofunction compactifications realise a generalised Satake compactification.  In the final section we briefly discuss this problem.

Horofunction compactifications of  Finsler symmetric spaces have been studied in a variety of  settings. Friedland and Freitas \cite{FF} determined the horofunction compactification for Finsler $p$-metrics on  $\mathrm{GL}_n(\mathbb{C})/\mathrm{U}_n$ for $1\leq p<\infty$.  The case $p=\infty$ on  $\mathrm{GL}_n(\mathbb{C})/\mathrm{U}_n$ corresponds to the Thompson distance on the symmetric cone $\Pi_n(\mathbb{C})$, and is an example of the spaces discussed here. For general normed spaces the horofunction compactification has been studied by Walsh in \cite{Wa2, Wa4}, who gave a complete description of the Busemann points. The horofunction for classical $\ell_p$-spaces have been analysed by Guti\`errez \cite{Gu2,Gu3}, see also \cite{FF}. For normed spaces with a polyhedral unit ball, the horofunction boundary has been determined in \cite{CKS,KMN,JS1}, and it was shown in \cite{JS} that they are related to projective toric varieties.  In \cite[Chapter 3]{Sc} Schilling showed that the horofunction compactification is homeomorphic to the closed dual unit ball for various classes of normed spaces including certain $l_1$-sums of normed spaces.  The Thompson distance horofunctions on general finite dimensional cones were studied by Walsh  in \cite{Wa3}, who gave a general charaterisation  of the Busemann points and the detour distance. In  Section \ref{thompson} we will give an explicit description of the Thompson distance horofunctions and the detour distance on symmetric cones in terms of the underlying Jordan algebra. As we shall see, the Jordan algebra structure and the symmetry of the space will allow us to give simple and direct proofs of these characterisations.  The horofunctions for general Hilbert distance spaces were  studied in \cite{Wa1}, see also \cite{Cl}. For Riemannian symmetric  spaces the horofunction compactification is homeomorphic to the Euclidean ball, see \cite{Eb}. In fact, for $\mathrm{CAT}(0)$ spaces it is known that the horofunction compactifications coincides with the visual boundary, see \cite{Ba,BH}.

\section{Preliminaries}
To fix the terminology and notation we recall the basic concepts and results concerning the horofunction compactification, the Hilbert and Thompson distances on cones,  and the theory of symmetric cones. 
\subsection{Horofunction compactifications} 
Let $(M,d)$ be a metric space and  $\mathbb{R}^M$ be the space of all real functions on $M$ equipped with the  topology of pointwise convergence. Fix  a  {\em basepoint} $b\in M$, and let $\mathrm{Lip}^c_b(M)$ denote the set of all $c$-Lipschitz functions $h\in\mathbb{R}^M$ such that $h(b)=0$, so $|h(x)-h(y)|\leq cd(x,y)$ for all $x,y\in M$. Then $\mathrm{Lip}^c_b(M)$ is a compact subset of $\mathbb{R}^M$. 
For $y\in M$ define the real valued function, 
\begin{equation}\label{internalpoint}
h_{y}(z) = d(z,y)-d(b,y)\mbox{\quad with $z\in M$.}
\end{equation}
Then $h_y(b)=0$ and $|h_y(z)-h_y(w)| = |d(z,y)-d(w,y)|\leq d(z,w)$. Thus, $h_y\in  \mathrm{Lip}_b^1(M)$ for all $y\in M$.  Using the fact that $\mathrm{Lip}_b^1(M)$ is compact, one defines the {\em horofunction compactification} of  $(M,d)$ to be the closure of $\{h_y\colon y\in M\}$ in $\mathbb{R}^M$, which we denote by $\overline{M}^h$. Its elements are called {\em metric functionals}, and the boundary $\partial \overline{M}^h= \overline{M}^h\setminus \{h_y\colon y\in M\}$ is called the {\em horofunction boundary}. The metric functionals in $\partial \overline{M}^h$ are called {\em horofunctions}. 

The topology of pointwise convergence on $ \mathrm{Lip}_b^1(M)$ coincides with the topology of uniform convergence on compact sets, see \cite[Section 46]{Mun}.  If $(M,d)$ is separable, then the topology of pointwise convergence on  $ \mathrm{Lip}_b^1(M)$ is  metrisable, and in that case each horofunction is the limit of a sequence $(x^n)$ in $M$.  If $(M,d)$ is proper (i.e. closed balls are compact) and geodesic, then the embedding $\iota\colon y\in M\mapsto h_y\in \mathbb{R}^M$ is a homeomorphism from $M$ onto $\iota(M)$, and $\overline{M}^h$ is a compactification in the usual topological sense.  Recall that a map $\gamma$ from a (possibly unbounded)  interval $I\subseteq \mathbb{R}$ into a metric space $(M,d)$ is called a {\em geodesic path} if 
\[
d(\gamma(s),\gamma(t)) = |s-t|\mbox{\quad for all }s,t\in I.
\]
The image, $\gamma(I)$, is called a {\em geodesic}, and a metric space $(M,d)$ is said to be {\em geodesic} if for each $x,y\in M$ there exists a geodesic path $\gamma\colon [a,b]\to M$ connecting $x$ and $y$, i.e, $\gamma(a)=x$ and $\gamma(b) = y$.  

The following well-known fact will be useful in the sequel. see \cite[Theorem 4.7]{Ri}. 
\begin{lemma}\label{Rieffel} 
If $(M,d)$ is a proper geodesic metric space, then $h\in  \partial \overline{M}^h$ if and only if there exists a sequence $(x^n)$ in $M$ with $d(b,x^n)\to \infty$ such that $(h_{x^n})$ converges to $h\in \overline{M}^h$ as $n\to \infty$.
\end{lemma}

A path $\gamma\colon T\to M$, where $T\subseteq [0,\infty)$ is unbounded and $0\in T$, is called an {\em almost geodesic} if  for all $\epsilon>0$ there exists an $M\geq 0$ such that 
\[
|d(\gamma(t),\gamma(s)) +d(\gamma(s),\gamma(0)) - t| <\epsilon\mbox{\quad for all $s,t\in T$ with $t\geq s\geq M$.}
\]
The notion of an almost geodesic sequence  goes back to  Rieffel \cite{Ri} and was further developed by Walsh and co-workers in \cite{AGW,LW,Wa4}. 
In particular, every almost geodesic  yields a horofunction for a proper geodesic metric space \cite[Lemma 4.5]{Ri}. 
\begin{lemma}\label{Rieffel2} Let $(M,d)$ be a proper geodesic metric space. If $\gamma$ is an almost geodesic in $M$, then 
\[
h(z) = \lim_{t\to\infty} d(z,\gamma(t))-d(b,\gamma(t))
\]
exists for all $z\in M$ and $h\in\partial \overline{M}^h$.
\end{lemma}
A horofunction $h$ in a  proper geodesic metric space $(M,d)$, is called a {\em Busemann point} if it is the limit of an almost geodesic $\gamma$ in $M$, and we denote the collection of all Busemann points by $\mathcal{B}_M$. 

Suppose that  $h$ and $h'$ are horofunctions of a proper geodesic metric space $(M,d)$. Let $W_h$ be the collection of neighbourhoods of $h$ in $\overline{M}^h$. The {\em detour cost} is defined by  
\[
H(h,h') = \sup_{W\in W_h}\left(\inf_{x\colon \iota(x)\in W} d(b,x)  +h'(x)\right),
\]
and the {\em detour distance} is given by $\delta(h,h') = H(h,h')+H(h',h)$. 

It is known, see for instance \cite[Lemma 3.1]{LW}, that  if $\gamma\colon T\to M$ is an almost geodesic converging to a horofunction $h$, then 
\begin{equation}\label{detourcost}
H(h,h') = \lim_{t\to\infty} d(b,\gamma(t)) +h'(\gamma(t))
\end{equation}
for all horofunctions $h'$. It is also  known, that on the set of Busemann points $\mathcal{B}_M$ the detour distance is a metric, where points can be at infinite distance from each other, see  e.g. \cite[Proposition 3.2]{LW}.  The detour distance yields a partition of $\mathcal{B}_M$ into equivalence classes,  called {\em parts}, where $h$ and $h'$ are equivalent if  $\delta(h,h') <\infty$. 

\subsection{Cones and the Hilbert and Thompson distances} 
Let $A$ be a real vector space. A {\em cone} $A_+$ in $A$ is a convex subset such that  $\lambda A_+\subseteq A_+$ for all $\lambda\geq 0$ and $A_+\cap -A_+=\{0\}$. A cone $A_+$ induces a partial ordering $\leq $ on $A$ by $x\leq y$ if $y-x\in A_+$. We write $x<y$ if $x\leq y$ and $x\neq y$.  The cone $A_+$ is said to be {\em Archimedean} if for each $x\in A$ and $y\in A_+$ with $nx\leq y$ for all $n\geq 1$, we have that $x\leq 0$. A point $u\in A_+$ is called an {\em order-unit} if for each $x\in A$ there exists $\lambda\geq 0$ such that $-\lambda u\leq x\leq \lambda u$. The triple $(A,A_+,u)$, where $A_+$ is an Archimedean cone and $u$ is an order-unit, is called an {\em order-unit space}, see \cite[Chapter 1]{AS0}. An order-unit space can be equipped with the {\em order-unit  norm}, 
\[
\|x\|_u =\inf\{\lambda\geq 0\colon -\lambda u\leq x\leq \lambda u\}.
\]
 So for each $x\in A$ we have that $-\|x\|_uu \leq x\leq \|x\|_uu$. Moreover, the cone $A_+$ is closed with respect to  the order-unit norm. The interior of $A_+$ will be denoted by $A_+^\circ$, and is nonempty, as  $u\in A_+^\circ$. 

Given $x\in A$ and $y\in A_+$, we say that $y$ {\em dominates} $x$ if there exist $\alpha,\beta\in\mathbb{R}$ such that $\alpha y\leq x\leq \beta y$. In that case, we let
\[
M(x/y) =\inf\{\beta\in\mathbb{R}\colon x\leq \beta y\}
\mbox{\quad and\quad }
m(x/y) = \sup\{\alpha\in\mathbb{R}\colon \alpha y\leq x\}.
\]

We note that for each  $x\in A$ we have that $\|x\|_u  = \max\{M(x/u), M(-x/u)\}$.
If $w\in A^\circ_+$, then $w$ dominates each $x\in A$, as it is an order-unit.  In that case we define 
\[|x|_w = M(x/w)- m(x/w)\mbox{\quad for } x\in A.\]
It can be shown that  $|\cdot|_w$ is a semi-norm on $A$, see \cite[Lemma A.1.1]{LNBook}, and  $|x|_w=0$ if and only if $x=\lambda w$ for some $\lambda\in\mathbb{R}$. So, $|\cdot|_w$ is a genuine norm on the quotient space $A/\mathbb{R}w$. Furthermore, 
\begin{equation}\label{seminorm}
|x+\lambda w|_w =|x|_w\mbox{\quad  for all $\lambda\in\mathbb{R}$ and $x\in V$. }
\end{equation}

The domination relation yields an equivalence relation on $A_+$ by $x\sim y$ if $y$ dominates $x$ and $x$ dominates $y$. So, $x\sim y$ if and only if there exist $0<\alpha\leq \beta$ such that $\alpha y\leq x\leq \beta y$. The equivalence classes are called  {\em parts} of $A_+$. The parts of a cone in a finite dimensional order-unit space correspond to the relative interiors of its faces, see \cite[Lemma 1.2.2]{LNBook}. 

The {\em Thompson distance} on $A_+$ is defined follows: 
\[
d_T(x,y)  =  \max\{\log M(x/y),\log M(y/x)\} 
\]
for all $x\sim y$ with $y\neq 0$, $d_T(0,0)=0$, and $d_T(x,y)=\infty$ otherwise.  The Thompson distance is a metric on each part of $A_+$, see \cite[Chapter 2]{LNBook}. In particular, it is a metric on $A_+^\circ$. 

Likewise, the {\em Hilbert distance}  on $A_+$ is given by, 
\[
d_H(x,y) =  \log M(x/y) +\log M(y/x) 
\]
for all $x\sim y$ with $y\neq 0$, $d_H(0,0)=0$, and $d_H(x,y)=\infty$ otherwise.  We note that $d_H(\lambda x,\mu y) = d_H(x,y)$ for all $x,y\in A_+$ and $\lambda,\mu>0$, so $d_H$ is not a metric on each part. It is, however, a metric between pairs of rays in each part of $A_+$, see \cite[Proposition 2.1.1]{LNBook}. 
 
In the sequel the following fact will be useful. The function $(x,y)\mapsto M(x/y)$ is continuous on $A_+\times A_+^\circ$, see \cite[Lemma 2.2]{LLNW}.

\subsection{Symmetric cones and Euclidean Jordan algebras}
The interior $A^\circ_+$ of a cone $A_+$ in a finite dimensional real vector space $A$ is called a {\em symmetric cone} if 
\begin{enumerate}
\item[(1)] there exists an inner product $(\cdot\mid\cdot)$ on $A$ such that $A_+$ is {\em self-dual}, i.e., \[A_+=A_+^*=\{y\in A\colon (y\mid x)\geq 0\mbox{ for all }x\in A_+\}.\]  
\item[(2)] $A^\circ_+$ is {\em homogeneous}, i.e.,  $\mathrm{Aut}(A_+)=\{g\in\mathrm{GL}(A)\colon g(A_+)=A_+\}$
acts transitively on $A_+^\circ$. 
\end{enumerate}  

The Koecher-Vinberg theorem say that there is a one-to-one correspondence between the symmetric cones  and the interiors of the cones of squares in finite dimensional formally real Jordan algebras (with unit). Recall that  a {\em real Jordan algebra} is a real vector space $A$ with a commutative bilinear  product $(x,y)\in A\times A\mapsto x\bullet y\in A$ satisfying the {\em Jordan identity}, $x^2\bullet (x\bullet y) = x\bullet (x^2\bullet y)$.
A Jordan algebra $A$  is said to be {\em formally real} if $x^2+y^2=0$ implies $x=0$ and $y=0$. 

Throughout this article we will  denote the unit in the Jordan algebra  by $u$. The unit $u$ is an order-unit for the cone $A_+=\{x^2\colon x\in A\}$.  We will consider the formally real  Jordan algebras $A$ as order-unit spaces $(A,A_+,u)$, where $A_+$ is the cone of squares and $u$ is the unit, and equip it with the order-unit norm. These normed spaces are precisely the finite dimensional  JB-algebras, see \cite[Theorem 1.11]{AS1}. 
 
Throughout we will fix the rank of the Jordan algebra $A$ to be $r$.  In a finite dimensional formally real Jordan algebra $A$ each $x\in A$ has a   {\em spectrum},  $\sigma(x)=\{ \lambda\in\mathbb{R}\colon \lambda u -x\mbox{ is not invertible}\}$, and we have that $A_+ =\{x\in V\colon \sigma(x)\subset [0,\infty)\}$. We write $\Lambda(x) =M(x/u)$ and note that 
\[\Lambda(x)=\max\{\lambda\colon \lambda\in\sigma(x)\},\]
 so that  
 \[
 \|x\|_u= \max\{\Lambda(x),\Lambda(-x)\} = \max\{|\lambda|\colon \lambda \in\sigma(x)\}
 \]
 for all $x\in A$. We also note that $\Lambda(x+\mu u) = \Lambda(x)+\mu$ 
 for all $x\in A$ and $\mu\in\mathbb{R}$. Moreover, if $x\leq y$, then $\Lambda(x)\leq\Lambda(y)$.

Recall that $p\in A$ is an {\em idempotent} if $p^2=p$. If, in addition, $p$ is non-zero and cannot be written as the sum of two non-zero idempotents, then it is said to be a {\em primitive} idempotent. The set of all primitive idempotent is denoted $\mathcal{J}_1(A)$ and is known to be a compact set \cite{Hi}. Two idempotents $p$ and $q$ are said to be orthogonal if $p\bullet q=0$, which is equivalent to $(p|q)=0$.  According to the spectral theorem \cite[Theorem III.1.2]{FK}, 
each $x$ has a {\em spectral decomposition}, $x = \sum_{i=1}^r \lambda_i p_i$, where each $p_i$ is a primitive idempotent, the $\lambda_i$'s are the eigenvalues of $x$ (including multiplicities), and $p_1,\ldots,p_r$ is a Jordan frame, i.e., the $p_i$'s are mutually orthogonal and $p_1+\cdots+p_r =u$. Moreover, $\mathrm{tr} x = \sum_{i=1}^r \lambda_i$ and $\det x = \prod_{i=1}^r \lambda_i$. 
 
For $x\in A$ we denote the {\em quadratic representation} by $U_x\colon A\to A$, which is the linear map, $U_x y = 2 x\bullet (x\bullet y) - x^2\bullet y$.  If $x\in A$ is invertible, then  $U_x$ is invertible and $U_x(A_+) = A_+$, see \cite[Proposition III.2.2]{FK}. 

Given a Jordan frame $p_1,\ldots, p_r$  in $A$ and $I\subseteq\{1,\ldots, r\}$ nonempty, we write $p_I = \sum_{i\in I} p_i$ and we let $A(p_I) = U_{p_I}(A)$.  Recall \cite[Theorem IV.1.1]{FK}  that $A(p_I)$ is the Peirce 1-space of the idempotent $p_I$, i.e, $A(p_I)=\{x\in A\colon p_I\bullet x =x\}$, which is a subalgebra. Given $z\in A(p_I)$, we write $\Lambda_{A(p_I)}(z)$ to denote the maximal eigenvalue of $z$ in the Jordan subalgebra $A(p_I)$ with unit $p_I$. So $ \Lambda_{A(p_I)}(z) =\inf\{\lambda\colon z\leq \lambda p_I\}$. 

\section{Thompson distance horofunctions}\label{thompson}
Let $A_+^\circ$ be a symmetric cone. For $x,y\in A_+^\circ$ we have that $x\leq \lambda y$ if and only if $U_{y^{-1/2}}x\leq \lambda u$. Therefore 
$\log M(x/y) = \log  \max \sigma(U_{y^{-1/2}}x) =  \max \sigma(\log U_{y^{-1/2}}x)$. We also have that  
\[
  \inf \{ \lambda > 0\colon y\leq \lambda x  \}
 =  (\sup \{ \mu > 0 \colon  \mu y \leq x \})^{-1}
 =(\sup \{ \mu > 0 \colon \mu u\leq  U_{y^{-1/2}}x \})^{-1}
= (\min \sigma(U_{y^{-1/2}}x))^{-1},
\]
so that $\log M(y/x) = \log (\min \sigma(U_{y^{-1/2}}x))^{-1} = - \min \sigma( \log U_{y^{-1/2}}x)$.
So the Thompson distance satisfies
\[
d_T(x,y) =\max\{  \max \sigma(\log U_{y^{-1/2}}x),  -\min \sigma(\log U_{y^{-1/2}}x)\} = \|\log U_{y^{-1/2}}x\|_u.
\]

The symmetry at $x\in A_+^\circ$ is given by $S_x(y) = U_xy^{-1}$ for $y\in A_+^\circ$. It can be shown that
\begin{equation}\label{inv}
M(x^{-1}/y^{-1}) = M(y/x)\mbox{\quad  for $x,y\in A^\circ_+$},
\end{equation}
see \cite[p.\,1518]{LRW}. Thus, for each $x\in A_+^\circ$ the symmetry $S_x\colon A_+^\circ\to A_+^\circ$ is a $d_T$-isometry.   

The following observation will be useful when determining the Thompson distance horofunctions for symmetric cones. 
\begin{lemma}\label{isom} Let $p_1,\ldots,p_k$ be orthogonal primitive idempotents in a finite dimensional formally real Jordan algebra $A$. 
The restriction of  $\mathrm{exp}_u\colon (A,\|\cdot\|_u)\to (A_+^\circ,d_T)$ to $\mathrm{Span}(\{p_1,\ldots,p_k\})$ is an isometry. Moreover, if $x,w\in \mathrm{Span}(\{p_1,\ldots,p_k\})$ with $\|w\|_u=1$ and $\gamma(t) = tw+x$ for $t\in\mathbb{R}$, then $\psi\colon t\mapsto\mathrm{exp}_u(\gamma(t))$ is a geodesic in $(A_+^\circ,d_T)$. 
\end{lemma}
\begin{proof}
For $y,z\in \mathrm{Span}(\{p_1,\ldots,p_k\})$ we have that 
\[
d_T(\mathrm{exp}_u(y),\mathrm{exp}_u(z))  =  d_T(e^y,e^z) =  \|\log U_{e^{-z/2}}e^y\|_u =\|\log e^{y-z}\|_u= \|y-z\|_u.
\]
The second assertion now easily follows, as $d_T(\mathrm{exp}_u(\gamma(t)),\mathrm{exp}_u(\gamma(s)))= \|(t-s)w\|_u = |t-s|$ for all $s,t\in\mathbb{R}$. 
\end{proof}

We now give a complete description of the Thompson distance horofunctions on symmetric cones. 
\begin{theorem}\label{horofunctions} If $A_+^\circ$ is a symmetric cone, then the  horofunctions of $(A_+^\circ,d_T)$ are all the functions of the form 
\begin{equation}\label{horo}
h(x) =  \max\{\log M(y/x),\log M(z/x^{-1})\}\mbox{\quad for }x\in A_+^\circ,
\end{equation}
where $y,z\in A_+$ with $y\bullet z =0$ and $\max \{\|y\|_u,\|z\|_u\}=1$. Here we use the convention that if $y$ or $z$ is $0$, then the corresponding term is omitted from the maximum. Moreover, each horofunction is a Busemann point.
\end{theorem}
\begin{proof}
We first show that each horofunction is of the from (\ref{horo}). Let $(y^n)\in A_+^\circ$ be such that $h_{y^n}\to h$ where $h$ is a horofunction. 
By Lemma \ref{Rieffel} we know that $d_T(y^n,u)\to\infty$. Let $r_n = e^{d_T(y^n,u)}$ and $z^n = (y^n)^{-1}$ for all $n$. Then $y^n\leq r_n u$ and $z^n\leq r_n u$.  Set $\hat{y}^n = y^n/r_n$ and $\hat{z}^n = z^n/r_n$. By taking a subsequence we may assume that $\hat{y}^n\to y$ and $\hat{z}^n\to z$, as $\hat{y}^n,\hat{z}^n\leq u$ for all $n$. As $\hat{y}^n\bullet \hat{z}^n = u/r_n^2 \to 0$, we conclude that $y\bullet z =0$. 

Also note that as $\|y^n\|_u = M(y^n/u)$ and $\|z^n\|_u = M(z^n/u) = M(u/y^n)$, we get that 
\[
r_n = e^{d_T(y^n,u)} = \max\{ M(y^n/u),M(u/y^n)\} = \max\{\|y^n\|_u,\|z^n\|_u\},
\]
so that $\max\{\|y\|_u,\|z\|_u\} = \max\{\|\hat{y}^n\|_u,\|\hat{z}^n\|_u\} = 1$. 
Using (\ref{inv}) we find  for $x\in A_+^\circ$ that 
\begin{eqnarray*}
h(x)  & = & \lim_{n\to\infty}  \max\{ \log M(y^n/x),\log M(x/y^n)\} - \log r_n\\
  & = & \lim_{n\to\infty} \max\{ \log (r_n^{-1}M(y^n/x)), \log (r_n^{-1}M(z^n/x^{-1}))\}\\
  & = & \lim_{n\to\infty} \max\{ \log M(\hat{y}^n/x), \log M(\hat{z}^n/x^{-1})\}. 
\end{eqnarray*}
Note that if $w^n,v\in A_+^\circ$ with $w^n\to 0$, then $M(w^n/v)\to M(0/v) = 0$ by the continuity of the $M$ function (\cite[Lemma 2.2]{LLNW}). 
As $h(x)\geq -d_T(x,u)$ for all $x\in A_+^\circ$, we deduce from the previous equality that 
$h(x)= \max\{ \log M(y/x), \log M(z/x^{-1})\}$, where if $y=0$ or $z=0$,  the corresponding term is omitted from the maximum.

To show that each function of the form (\ref{horo}) is a horofunction, let $y,z\in A_+$ with $\max\{\|y\|_u,\|z\|_u\}=1$ and $y\bullet z =0$. We will discuss the case where $y$ and $z$ are both non-zero. The other cases can be proved in the same fashion and  are left  to the reader.

Using the spectral decomposition we can write 
\[
y = \sum_{i\in I} e^{-\alpha_i}p_i \mbox{\quad and \quad } z= \sum_{j\in J} e^{-\alpha_j}p_j,
\]
where $\min\{\alpha_k\colon k\in I\cup J\}=0$, and the $p_k$'s are mutually orthogonal primitive idempotents. Set $p_I = \sum_{i\in I} p_I$ and $p_J = \sum_{j\in J} p_j$. So $p_I\bullet p_J =0$, and hence orthogonal idempotents in $A$.  
Let $v = -\sum_{i\in I} \alpha_i p_i +\sum_{j\in J}\alpha_j p_j$ and $w = p_I-p_J$. Note that $w\in \mathrm{Span}(\{p_k\colon k\in I\cup J\})$. For $t>0$ let $\gamma(t) = tw+v$. 

From Lemma \ref{isom} we know that $\psi\colon t\mapsto \mathrm{exp}_u(\gamma(t))$ is a geodesic in $A_+^\circ$, and  for all $t>0$ sufficiently large we have that 
\[
d_T(\psi(t),u) = \|\log e^{tw+v}\|_u = \|tw+v\|_u = \max\{|t-\alpha_k|\colon k\in I\cup J\} = t,
\]
as $\min\{\alpha_k\colon k\in I\cup J\}=0$. Moreover, $\psi(t)^{-1} = e^{-(tw+v)}$. 

Thus, 
\[
\lim_{t\to\infty} e^{-t}\psi(t) = \lim_{t\to\infty} \sum_{i\in I} e^{-\alpha_i}p_i +\sum_{j\in J} e^{-2t +\alpha_j}p_j + e^{-t}(u-p_I-p_J)= y
\]  
and 
\[
\lim_{t\to\infty} e^{-t}\psi(t)^{-1} = \lim_{t\to\infty} \sum_{i\in I} e^{-2t+\alpha_i}p_i +\sum_{j\in J} e^{-\alpha_j}p_j + e^{-t}(u-p_I-p_J)= z.
\]  
Now using (\ref{inv}) we deduce  for each $x\in A_+^\circ$ that 
\begin{eqnarray*}
\lim_{t\to\infty} d_T(x,\psi(t))-d_T(\psi(t),u) & = & \lim_{t\to\infty} \max\{\log M(\psi(t)/x), \log M(x/\psi(t))\} - t\\
  & = & \lim_{t\to\infty} \max\{\log M(e^{-t}\psi(t)/x), \log M(e^{-t}\psi(t)^{-1}/x^{-1})\} \\
  & = & \max\{\log M(y/x), \log M(z/x)\}, 
\end{eqnarray*}
as $v\in A_+\mapsto M(v/w)$ is continuous for all $w\in A_+^\circ$ by \cite[Lemma 2.2]{LLNW}.  We conclude that the function of the form (\ref{horo}) is a horofunction. In fact, it is a Busemann point, as $t\mapsto \psi(t)$ is a geodesic. 
\end{proof}
In the case of $\Pi_n(\mathbb{C})$ we get that the horofunctions are the functions $h\colon \Pi_n(\mathbb{C})\to\mathbb{R}$ of the form: 
\[
h(X) = \max\{\log \max\sigma(X^{-1/2}AX^{-1/2}),\log \max\sigma(X^{1/2}BX^{1/2})\},
\]
where $A$ and $B$ are positive semi-definite, with $\max\{\max\sigma(A), \max\sigma(B)\}=1$ and $1/2(AB+BA)=0$. 

\begin{remark}
We see from the proof of Theorem \ref{horofunctions}  that each horofunction of $(A_+^\circ,d_T)$  is obtained as the limit of a geodesic in the span of a Jordan frame.  This is in agreement with  \cite[Lemma 4.4]{HSWW}, which implies that each horofunction of $(A_+^\circ,d_T)$  arises as a limit of a sequence in a flat.   In fact, if $h$ is a horofunction given by $h(x) =  \max\{\log M(y/x),\log M(z/x^{-1})\}$, where $ y =\sum_{i\in I}e^{-\alpha_i}p_i$, $ z=\sum_{j\in J}e^{-\alpha_j}p_j$,  then $\psi\colon t\mapsto \mathrm{exp}_u(\gamma(t))$, with 
\[
\gamma(t) = t\left(\sum_{i\in I} p_i - \sum_{j\in J} p_j\right) - \sum_{i\in I} \alpha_ip_i +\sum_{j\in J}\alpha_jp_j
\]
for $t>0$, is a geodesic $A_+^\circ$ converging to $h$. 
\end{remark}

To analyse the parts and the detour distance,  the following observation is useful. If $(V_i,C_i,u_i)$, $i=1,2$, are order-unit spaces, then the product space $V_1\oplus V_2$ is an order unit space with cone $C_1\times C_2$ and order unit $u=(u_1,u_2)$. Moreover, for   $x=(x_1,x_2), y=(y_1,y_2)\in C_1\times C_2$ we have that 
\[
M(x/y) = \max\{M(x_1/y_1),M(x_2/y_2)\}.
\]
Indeed, if $x\leq \lambda y$, then $\lambda y -x \in C_1\times C_2$, and hence $\lambda y_1-x_1\in C_1$ and $\lambda y_2-x_2\in C_2$. So $\lambda\geq \max\{M(x_1/y_1),M(x_2/y_2)\}$. On the hand, if $x_1\leq \mu y_1$ and $x_2\leq \mu y_2$, then $\mu y-x\in C_1\times C_2$, so that 
$\mu\geq M(x/y)$. 
\begin{theorem}\label{detour}
Suppose that  $h$ and $h'$ are horofunctions in $\overline{A_+^\circ}^h$, where
\[
h(x) = \max\{\log M(y/x),\log M(z/x^{-1})\} \mbox{ and } 
h'(x) =  \max\{\log M(y'/x),\log M(z'/x^{-1})\}\mbox{\quad for }x\in A_+^\circ.
\] 
Let $y=\sum_{i\in I} e^{-\alpha_i}p_i$, $z =\sum_{j\in J} e^{-\alpha_j}p_j$, and set $p_I=\sum_{i\in I} p_i$ if $y\neq 0$, and $p_J=\sum_{j\in J} p_j$ if $z\neq 0$. 
Then $h$ and $h'$ are in the same part if and only if $y\sim y'$ and $z\sim z'$.  Moreover, in that case, 
\[
\delta(h,h') = d_H((y,z),(y',z'))\mbox{,\quad where  $(y,z), (y',z')\in U_{p_I}(A)\oplus U_{p_J}(A)$}
\]
and $d_H$ is the Hilbert distance on the product cone $U_{p_I}(A)_+\times U_{p_J}(A)_+$.  Here, if $y$ or $z$ is $0$, we omit the corresponding term in the sum $U_{p_I}(A)\oplus U_{p_J}(A)$. 
\end{theorem}
\begin{proof}
We use (\ref{detourcost}) to determine the detour cost $H(h,h')$.   Note that as $\|y\|_u=1$ or $\|z\|_u=1$, we have that $\min\{\alpha_k\colon k\in I\cup J\}=0$.  Using the notation as in the proof of Theorem \ref{horofunctions} we let 
$v = -\sum_{i\in I} \alpha_i p_i +\sum_{j\in J}\alpha_j p_j$ and $w = p_I-p_J$, so  $w\in  \mathrm{Span}(\{p_k\colon k\in I\cup J\})$. For $t>0$ let $\gamma(t) = tw+v$. 
From Lemma \ref{isom} we know that $\psi\colon t\mapsto \mathrm{exp}_u(\gamma(t))$ is a geodesic in $A_+^\circ$, and  for all $t>0$ sufficiently large we have that 
\[
d_T(\psi(t),u) = \|\log e^{tw+v}\|_u = \|tw+v\|_u = \max\{|t-\alpha_k|\colon k\in I\cup J\} = t,
\]
as $\min\{\alpha_k\colon k\in I\cup J\}=0$. Moreover, $\psi(t)^{-1} = e^{-(tw+v)}$ and $h_{\psi(t)}\to h$. 

As $e^{-t}\psi(t) = \sum_{i\in I} e^{-\alpha_i}p_i +\sum_{j\in J} e^{-2t +\alpha_j}p_j + e^{-t}(u-p_I-p_J)$, we have that $e^{-t}\psi(t)\to y$ and 
\[
e^{-t}\psi(t) \leq e^{-s}\psi(s)\mbox{\quad for all } 0\leq s\leq t. 
\]
Likewise, $e^{-t}\psi(t)^{-1}\to z$ and $e^{-t}\psi(t)^{-1} \leq e^{-s}\psi(s)^{-1}$ for all $0\leq s\leq t$. 

It follows from \cite[Lemma 5.5]{LP} that  if $y$ dominates $y'$ and $z$ dominates $z'$, then 
\begin{eqnarray*}
H(h,h') & = & \lim_{t\to\infty} d_T(\psi(t),u) + h'(\psi(t)) \\
   & = & \lim_{t\to\infty} t + \max \{\log M(y'/\psi(t)),\log M(z'/\psi(t)^{-1}) \}\\
   & = & \lim_{t\to\infty} \max \{\log M(y'/e^{-t}\psi(t)),\log M(z'/e^{-t}\psi(t)^{-1}) \}\\
   & = & \max \{\log M(y'/y),\log M(z'/z)\},
\end{eqnarray*}
and otherwise, $H(h,h')=\infty$. Here, if $y'=0$ or $z'=0$, the corresponding term is omitted from the maximum.

Interchanging the roles of $h$ and $h'$ gives 
\[
\delta(h,h') =  \max\{\log M(y'/y),\log M(z'/z)\} +  \max\{\log M(y/y'),\log M(z/z')\} = d_H((y,z),(y',z'))
\] 
if $y\sim y'$ and $z\sim z'$, and $\delta(h,h')=\infty$ otherwise. 
\end{proof}

\section{Proof of part (a) of Theorem \ref{main}}
To define the extension of the exponential map we need to recall the description of the horofunctions in $\overline{A}^h$ from \cite[Theorem 4.2]{LP}. It was shown there that the horofunctions in $\overline{A}^h$ are precisely  the functions $g\colon A\to \mathbb{R}$ of the form,
\begin{equation}\label{horA}
g(v) = \max\Big\{\Lambda_{A(p_I)}(-U_{p_I}v -\sum_{i\in I}\alpha_ip_i),\Lambda_{A(p_J)}(U_{p_J}v -\sum_{j\in J}\alpha_jp_j)\Big\},
\end{equation}
where $p_1,\ldots,p_r \in A$ is a Jordan frame, $I,J\subseteq\{1,\ldots,r\}$ are disjoint and $I\cup J$ is nonempty, $p_I=\sum_{i\in I} p_i$, $p_J=\sum_{j\in J} p_j$, and $\alpha \in\mathbb{R}^{I\cup J}$ with $\min\{\alpha_k\colon k\in I\cup J\}$. Here the convention is that if $I$ or $J$ is empty, the corresponding term in the maximum is omitted.

\begin{definition} The {\em exponential map} $\mathrm{exp}_u\colon \overline{A}^h\to \overline{A_+^\circ}^h$ is defined by, $\mathrm{exp}_u(v) = e^v$ for $v\in A$, and for $g\in\partial \overline{A}^h$ given by (\ref{horA}) we let  $\mathrm{exp}_u(g) = h$, where 
\begin{equation}\label{Expg}
h(x)= \max\{\log M(y/z),\log M(z/x^{-1})\}\mbox{\quad for }x\in A^\circ_+,
\end{equation}
with $y= \sum_{i\in I}e^{-\alpha_i}p_i$ and $z=\sum_{j\in J}e^{-\alpha_j}p_j$. 
\end{definition}
Note that $y,z\in A_+$, with $\max\{\|y\|_u,\|z\|_u\}=1$, as $\min\{\alpha_k\colon k\in I\cup J\}=0$, and $y\bullet z=0$. So, $\mathrm{exp}_u(g)$ is a horofunction by Theorem \ref{horofunctions}. Moreover, the extension is well defined. To show this we use the following observation. 
\begin{lemma}\label{equal} Suppose $x,y\in A$ have spectral decompositions $x=\sum_{i\in I} \alpha_i p_i$ and $y = \sum_{j\in J} \beta_j q_j$. If $\sum_{i\in I} p_i =\sum_{j\in J}q_j$ and $x=y$, then $\sum_{i\in I} e^{-\alpha_i}p_i = \sum_{j\in J} e^{-\beta_j}q_j$. 
\end{lemma} 
\begin{proof} 
Let $p_I = \sum_{i\in I} p_i$ and $q_J = \sum_{j\in J} q_j$, so $p_I =q_J$. As   $\sum_{i\in I} \alpha_i p_i  =\sum_{j\in J} \beta_j q_j$,  we have that $(u-p_I) +\sum_{i\in I} e^{-\alpha_i} p_i  = (u-q_J)+\sum_{j\in J} e^{-\beta_j} q_j$, and hence $\sum_{i\in I} e^{-\alpha_i}p_i = \sum_{j\in J} e^{-\beta_j}q_j$. 
\end{proof}
 Now to see that the extension is well-defined assume that  $g$ in (\ref{horA}) is represented differently as 
\[
g(v) = \max\Big\{\Lambda_{A(q_{I'})}(-U_{q_{I'}}v -\sum_{i\in I'}\beta_iq_i),\Lambda_{A(q_{J'})}(U_{q_{J'}}v -\sum_{j\in J'}\beta_jq_j)\Big\},
\]
It follows from  \cite[Theorem 4.3]{LP}  and the fact that $\delta(g,g)=0$ that $p_I = q_{I'}$, $p_J=q_{J'}$, and 
 \[
 \sum_{i\in I} \alpha_ip_i +\sum_{j\in J}\alpha_jp_j = \sum_{i\in I'} \beta_iq_i +\sum_{j\in J'}\beta_jq_j, 
 \]
 as $\min\{\alpha_m\colon m\in I\cup J\} =0 =\min\{\beta_m\colon m\in I'\cup J'\}$. So, 
 \[
  \sum_{i\in I} \alpha_ip_i= U_{p_I}( \sum_{i\in I} \alpha_ip_i +\sum_{j\in J}\alpha_jp_j ) = U_{q_{I'}}(\sum_{i\in I'} \beta_iq_i +\sum_{j\in J'}\beta_jq_j) = \sum_{i\in I'} \beta_iq_i 
 \]
 and 
 \[
  \sum_{j\in J} \alpha_jp_j= U_{p_J}( \sum_{j\in J} \alpha_ip_i +\sum_{j\in J}\alpha_jp_j ) = U_{q_{J'}}(\sum_{i\in I'} \beta_iq_i +\sum_{j\in J'}\beta_jq_j) = \sum_{j\in J'} \beta_jq_j. 
 \]
 Using Lemma \ref{equal} we conclude that $\sum_{i\in I} e^{-\alpha_i}p_i = \sum_{i\in I'} e^{-\beta_i}q_i$  and $\sum_{j\in J}e^{-\alpha_j}p_j = \sum_{j\in J'}e^{-\beta_j}q_j$, 
 and hence the extension is well defined. 

We will prove in this section that this extension of the exponential map is a homeomorphism that maps each part in $\partial \overline{A}^h$ onto a part in $\partial \overline{A_+^\circ}^h$. To start we show that the extension is a bijection. 
\begin{lemma}\label{bijection1} The map $\mathrm{exp}_u\colon \overline{A}^h\to \overline{A_+^\circ}^h$ is a bijection, which maps $A$ onto $A_+^\circ$ and $\partial \overline{A}^h$ onto $\partial \overline{A_+^\circ}^h$.
\end{lemma}
\begin{proof}
Clearly  $\mathrm{exp}_u$ is a bijection from $A$ onto $A_+^\circ$. Moreover, it follows from Theorem \ref{horofunctions} and \cite[Theorem 4.2]{LP} that $\mathrm{exp}_u$ maps $\partial \overline{A}^h$ onto $\partial \overline{A_+^\circ}$.  To show that the extension is injective on $\partial \overline{A}^h$ suppose that $h=\mathrm{exp}_u(g) = \mathrm{exp}_u(g')=h'$, where $g$ is as in (\ref{horA}) and 
\begin{equation}\label{horB}
g'(v) = \max\Big\{\Lambda_{A(q_I{'})}(-U_{q_{I'}}v -\sum_{i\in I'}\beta_iq_i),\Lambda_{A(q_{J'})}(U_{q_{J'}}v -\sum_{j\in J'}\beta_jq_j)\Big\},
\end{equation}
with $q_1,\ldots,q_r\in A$ is a Jordan frame, $I',J'\subseteq\{1,\ldots,r\}$ are disjoint and $I'\cup J'$ is nonempty, and $\beta \in\mathbb{R}^{I\cup J}$ with $\min\{\beta_k\colon k\in I'\cup J'\}$. 
 
By definition of the extension, we have that 
\[h(x)= \max\{\log M(y/x),\log M(z/x^{-1})\},\mbox{\quad with $y= \sum_{i\in I}e^{-\alpha_i}p_i$ and $z=\sum_{j\in J}e^{-\alpha_j}p_j$},
\] 
and 
\[h'(x)= \max\{\log M(y'/x),\log M(z'/x^{-1})\},\mbox{\quad with $y'= \sum_{i\in I'}e^{-\beta_i}q_i$ and $z'=\sum_{j\in J'}e^{-\beta_j}q_j$. }
\] 

As $h=h'$ we know that $\delta(h,h')=0$, and hence $y=y'$ and $z=z'$ by Theorem \ref{detour}. It now follows from the spectral theorem \cite[Theorems III.1.1 and III.1.2]{FK}  that $ \sum_{i\in I}\alpha_ip_i=\sum_{i\in I'}\beta_iq_i$ and  $\sum_{j\in J}\alpha_jp_j=\sum_{j\in J'}\beta_jq_j$. 
Moreover, $p_I=q_{I'}$ and $p_J=q_{J'}$. This implies that $g=g'$, and hence $\mathrm{exp}_u$ is injective, which completes the proof.
\end{proof}

Clearly $\mathrm{exp}_u$ is continuous on $A$. To establish the continuity on all of $\overline{A}^h$ we prove two lemmas. 
\begin{lemma}\label{cont1} If $(w^n)$ in $A$ converges to $g\in\partial \overline{A}^h$, then $(\mathrm{exp}_u(w^n))$ converges to $\mathrm{exp}_u(g)$. 
\end{lemma}
\begin{proof}
To prove the statement we show that each subsequence of $(\mathrm{exp}_u(w^n))$  has a convergent subsequence with limit $\mathrm{exp}_u(g)$. 
So let $(\mathrm{exp}_u(w^{n_k}))$ be a subsequence and let $g$ be given by (\ref{horA}), so that $h=\mathrm{exp}_u(g)$ is given by (\ref{Expg}). As $g$ is a horofunction, we know by Lemma \ref{Rieffel} that $\|w^{n_k}\|_u\to\infty$. 
For $k\geq 1$ write $r_{n_k} =\|w^{n_k}\|_u$ and let $w^{n_k} = \sum_{i=1}^r \lambda^{n_k}_iq_i^{n_k}$ be the spectral decomposition of $w^{n_k}$. After taking a subsequence, may assume that  
\begin{enumerate}[(1)]
\item There exists $s\in\{1,\ldots,r\}$ such that $r_{n_k} = |\lambda^{n_k}_s|$ for all $k\geq 1$.
\item There exist $I_+\subseteq\{1,\ldots,r\}$ such that for each $k\geq 1$ we have $\lambda^{n_k}_i>0$ if and only if $i\in I_+$. 
\item $q^{n_k}_i\to q_i$ for all $i\in\{1,\ldots,r\}$.
\end{enumerate} 
The third property follows from the fact that the set of primitive idempotents is compact, see \cite{Hi}. 

Now let $\beta_i^{n_k} = r_{n_k}-\lambda^{n_k}_i$ for all $i\in I_+$ and $\beta_i^{n_k} = r_{n_k}+\lambda^{n_k}_i$ for all $i\not\in I_+$. 
So $\beta^{n_k}_i\geq 0$ for all $i$, and $\beta^{n_k}_s =0$ for all $k$. By taking a further subsequence we may also assume that 
$\beta^{n_k}_i\to \beta_i\in [0,\infty]$ for all $i$.  Let $I'=\{i\in I_+\colon \beta_i<\infty\}$ and $J'=\{j\not\in I_+\colon \beta_j<\infty\}$. Note that $s\in I'\cup J'$ and hence the union is nonempty. Moreover, $I'$ and $J'$ are disjoint.   
It now follows from \cite[Lemma 4.7]{LP} that $g_{w^{n_k}}\to g'$, where 
\[
g'(v) = \max\Big\{\Lambda_{A(q_{I'})}(-U_{q_{I'}}v -\sum_{i\in I'}\beta_iq_i),\Lambda_{A(q_{J'})}(U_{q_{J'}}v -\sum_{j\in J'}\beta_jq_j)\Big\},
\]
with $q_{I'} =\sum_{i\in I'} q_i$ and $q_{J'} = \sum_{j\in J'} q_j$. 

As $(w^n)$ converges to $g$, we find that $g=g'$, and hence $\delta(g,g')=0$. It now follows from \cite[Theorem 4.3]{LP} that 
$p_I = q_{I'}$, $p_J=q_{J'}$, and  
\[
\sum_{i\in I} \alpha_ip_i +\sum_{j\in J}\alpha_jp_j= \sum_{i\in I'} \beta_i q_i+\sum_{j\in J'} \beta_j q_j, 
\]
as $\min_{i\in I\cup J} \alpha_i =0= \min_{i\in I'\cup J'} \beta_i$
This implies that 
\[
\sum_{i\in I} \alpha_ip_i  = U_{p_{I}} (\sum_{i\in I} \alpha_ip_i +\sum_{j\in J}\alpha_jp_j) = U_{q_{I'}} ( \sum_{i\in I'} \beta_i q_i+\sum_{j\in J'} \beta_j q_j) =  \sum_{i\in I'} \beta_i q_i.
\]
Likewise, we have 
$\sum_{j\in J}\alpha_jp_j= \sum_{j\in J'} \beta_j q_j$. 
From Lemma \ref{equal}  we now deduce that $ y=\sum_{i\in I} e^{-\alpha_i}p_i  = \sum_{i\in I'} e^{-\beta_i} q_i$ and $z= \sum_{j\in J} e^{-\alpha_j}p_j  = \sum_{j\in j'} e^{-\beta_j} q_j$. 

Note that  
  \[ 
\lim_{k\to\infty} e^{-r_{n_k}} \mathrm{exp}_u(w^{n_k}) =  \lim_{k\to\infty} \sum_{i=1}^r e^{-(r_{n_k} -\lambda_i^{n_k})}q^{n_k}_i  = \sum_{i\in I'} e^{-\beta_i}q_i =y
\]
and 
 \[ 
\lim_{k\to\infty} e^{-r_{n_k}} \mathrm{exp}_u(-w^{n_k}) =  \lim_{k\to\infty} \sum_{i=1}^r e^{-(r_{n_k} +\lambda_i^{n_k})}q^{n_k}_i  = \sum_{j\in J'} e^{-\beta_j}q_j =z.
\]
It now follows from the continuity of the $M$ function, see \cite[Lemma 2.2]{LLNW}, that 
\begin{eqnarray*}
\lim_{k\to\infty} h_{\mathrm{exp}_u(w^{n_k})}(x) & = & \lim_{k\to\infty} d_T(x,\mathrm{exp}_u(w^{n_k})) - d_T(u,\mathrm{exp}_u(w^{n_k})) \\
   & = & \lim_{k\to\infty}\max\{\log M(\mathrm{exp}_u(w^{n_k})/x),\log M(\mathrm{exp}_u(-w^{n_k})/x^{-1})\} - \log e^{r_{n_k}}\\
   & = &\lim_{k\to\infty} \max\{\log M(e^{-r_{n_k}}\mathrm{exp}_u(w^{n_k})/x),\log M(e^{-r_{n_k}}\mathrm{exp}_u(-w^{n_k})/x^{-1})\}\\
   & = & \max\{\log M(y/x),\log M(z/x^{-1})\},
  \end{eqnarray*}
  which shows that $(\mathrm{exp}_u(w^{n_k}))$ converges to $h = \mathrm{exp}_u(g)$, and hence the proof is complete. 
 \end{proof}
 Next we show continuity  of $\mathrm{exp}_u$ in the  horofunction boundary.
 \begin{lemma}\label{cont2} If $(g_n)$ in $\partial \overline{A}^h$ converges to a horofunction $g$, then 
 $(\mathrm{exp}_u(g_n))$  converges to $\mathrm{exp}_u(g)$. 
 \end{lemma}
 \begin{proof}
 Let $(g_n)$ be a sequence in $\partial \overline{A}^h$ converging to $g$, where $g$ is given by (\ref{horA}). So $\mathrm{exp}_u(g)=h$, where is $h$ given by (\ref{Expg}). To show the lemma, we prove that each subsequence of $(\mathrm{exp}_u(g_n))$ has a convergent subsequence with limit $h$.  Let $(\mathrm{exp}_u(g_{n_k}))$ be a subsequence. 
 By \cite[Theorem 4.2]{LP} we can write  for $k\geq 1$, 
 \[
 g_{n_k}(v)  = \max\Big\{\Lambda_{A(q^k_{I_k})}(-U_{q^k_{I_k}}v -\sum_{i\in I_k}\beta^k_iq^k_i),\Lambda_{A(q^k_{J_k})}(U_{q^k_{J_k}}v -\sum_{j\in J_k}\beta^k_jq^k_j)\Big\},
 \]
 where the $q^k_i$ and $q^k_j$ are orthogonal primitive idempotents, $I_k,J_k\subseteq\{1,\ldots,r\}$ are disjoint with $I_k\cup J_k$ nonempty, 
 $\min\{\beta^k_m\colon m\in I_k\cup J_k\} =0$, $q_{I_k}^k = \sum_{i\in I_k} q_i^k$ and $q_{J_k}^k = \sum_{j\in J_k} q_j^k$. 
 
 The approach will be similar to the one taken in the proof of the previous lemma. After taking subsequences we may assume that:
 \begin{enumerate}[(1)]
 \item There exist $I_0,J_0\subseteq \{1,\ldots,r\}$ such that $I_0=I_k$ and $J_0=J_k$ for all $k$.
 \item There exists $s\in I_0\cup J_0$ such that $\beta_s^k=0$ for all $k$. 
 \item $\beta_m^k\to\beta_m\in [0,\infty]$ and $q^k_m\to q_m$ for all $m\in I_0\cup J_0$. 
 \end{enumerate}
 Now let $I'=\{i\in I_0\colon \beta_i<\infty\}$ and $J'=\{j\in J_0\colon \beta_j<\infty\}$, and note that $s\in I'\cup J'$.
 
Next we show that 
\begin{equation}\label{conv}
\lim_{k\to\infty} g_{n_k}(v) = \max\Big\{\Lambda_{A(q_{I'})}(-U_{q_{I'}}v -\sum_{i\in I'}\beta_iq_i),\Lambda_{A(q_{J'})}(U_{q_{J'}}v -\sum_{j\in J'}\beta_jq_j)\Big\}
 \end{equation}
 where the term is omitted if the corresponding set $I'$ or $J'$ is empty. Here $q_{I'} = \sum_{i\in I'} q_i$ and $q_{J'} = \sum_{j\in J'} q_j$.
 
 First let us assume that $I_0\neq \emptyset$ and $I'=\emptyset$. Note that 
 \[
 -U_{q^k_{I_0}}v\leq \|v\|_u U_{q^k_{I_0}}u=\|v\|_uq^k_{I_0},
 \] 
 as $-v\leq \|v\|_uu$ and $U_{q^k_{I_0}}(A_+)\subseteq A_+$. 
 It follows that 
 \[
 -U_{q^k_{I_0}}v - \sum_{i\in I_0} \beta_i^kq_i^k \leq \sum_{i\in I_0} (\|v\|_u -\beta_i^k)q_i^k,
 \]
 so that 
 \[
 \Lambda_{A(q^k_{I_0})}(-U_{q^k_{I_0}}v -\sum_{i\in I_0}\beta^k_iq^k_i)\leq \Lambda_{A(q^k_{I_0})}(\sum_{i\in I_0} (\|v\|_u -\beta_i^k)q_i^k)
 \leq \max_{i\in I_0} (\|v\|_u - \beta^k_i).
\]
The right-hand side diverges to $-\infty$ as $k\to\infty$, since $I'$ is empty.  We also know that for each horofunction $\bar{g}$ in $\overline{A}^h$ we have that $\bar{g}(v)\geq -\|v\|_u$.  So, if $I_0 \neq \emptyset$ and $I' =\emptyset$, then $s\in J'$ and for each $v\in A$ we have that 
\[
g_{n_k}(v) = \Lambda_{A(q^k_{J_0})}(U_{q^k_{J_0}}v -\sum_{j\in J_0}\beta^k_jq^k_j)
\]
for all $k$ large. In the same way we get that if $J_0\neq \emptyset$ and $J' =\emptyset$, then 
\[
g_{n_k}(v) = \Lambda_{A(q^k_{I_0})}(-U_{q^k_{I_0}}v -\sum_{i\in I_0}\beta^k_iq^k_i)
\]
for all $k$ large. 

On the other hand, if $I' \neq \emptyset$, then by \cite[Lemma 4.7]{LP} we know that 
\[
\Lambda_{A(q^k_{I_0})}(-U_{q^k_{I_0}}v -\sum_{i\in I_0}\beta^k_iq^k_i)\to \Lambda_{A(q_{I'})}(-U_{q_{I'}}v -\sum_{i\in I'}\beta_iq_i),  
\]
and, similarly,  if $J'\neq \emptyset$, we have 
\[
\Lambda_{A(q^k_{J_0})}(U_{q^k_{J_0}}v -\sum_{j\in J_0}\beta^k_jq^k_j)
\to
\Lambda_{A(q_{J'})}(U_{q_{J'}}v -\sum_{j\in J'}\beta_jq_j).
\]
Note that if $I_0=\emptyset$, then $s\in J'$. Likewise, if $J_0=\emptyset$, then $s\in I'$. It follows that in each of the cases (\ref{conv}) holds.

 Let $g'\colon A\to\mathbb{R}$ be given by 
 \[
 g'(v) =  \max\Big\{\Lambda_{A(q_{I'})}(-U_{q_{I'}}v -\sum_{i\in I'}\beta_iq_i),\Lambda_{A(q_{J'})}(U_{q_{J'}}v -\sum_{j\in J'}\beta_jq_j)\Big\}.
 \]
 Then by \cite[Theorem 4.2]{LP} we know that $g'\in \partial \overline{A}^h$. As $g_n\to g$, we find that $g=g'$ and hence $\delta(g,g')=0$. 
 
 It now follows from \cite[Theorem 4.3]{LP} that  $p_I = q_{I'}$, $p_J=q_{J'}$, and 
 \[
 \sum_{i\in I} \alpha_ip_i +\sum_{j\in J}\alpha_jp_j = \sum_{i\in I'} \beta_iq_i +\sum_{j\in J'}\beta_jq_j, 
 \]
 as $\min\{\alpha_m\colon m\in I\cup J\} =0 =\min\{\beta_m\colon m\in I'\cup J'\}$.  This implies that  
 \[
  \sum_{i\in I} \alpha_ip_i= U_{p_I}( \sum_{i\in I} \alpha_ip_i +\sum_{j\in J}\alpha_jp_j ) = U_{q_{i'}}(\sum_{i\in I'} \beta_iq_i +\sum_{j\in J'}\beta_jq_j) = \sum_{i\in I'} \beta_iq_i 
 \]
 and 
 \[
  \sum_{j\in J} \alpha_jp_j= U_{p_J}( \sum_{j\in J} \alpha_ip_i +\sum_{j\in J}\alpha_jp_j ) = U_{q_{J'}}(\sum_{i\in I'} \beta_iq_i +\sum_{j\in J'}\beta_jq_j) = \sum_{j\in J'} \beta_jq_j. 
 \]
 Using Lemma \ref{equal} we conclude that 
  \[
 \sum_{i\in I} e^{-\alpha_i}p_i = \sum_{i\in I'} e^{-\beta_i}q_i \mbox{\quad  and \quad }\sum_{j\in J}e^{-\alpha_j}p_j = \sum_{j\in J'}e^{-\beta_j}q_j,
 \]

 So, if we let $\bar{y}^{k} = \sum_{i\in I_0}e^{-\beta^k_i}q^k_i$ and $\bar{z}^k = \sum_{j\in J_0}e^{-\beta^k_j}q^k_j$, then 
 \[
 \lim_{k\to\infty} \bar{y}^{k}  =  \sum_{i\in I'}e^{-\beta_i}q_i =  y\mbox{\quad and\quad } 
 \lim_{k\to\infty} \bar{z}^{k}  =  \sum_{j\in J'}e^{-\beta_j}q_j =  z.
 \]
Using the continuity of the $M$ function (\cite[Lemma 2.2]{LLNW}) we now get that 
 \[
 \lim_{k\to\infty} \mathrm{exp}_u(g_{n_k})(x) = \lim_{k\to\infty} \max\{\log M(\bar{y}^k/x),\log M(\bar{z}^k/x^{-1}) \} = 
\max\{\log M(y/x),\log M(z/x^{-1}) \} = h(x),
 \]
 which completes the proof.
  \end{proof}
 To complete the proof Theorem \ref{main}(a) the following concepts are useful.  
For $x,z\in A$ we let $[x,z]=\{y\in A\colon x\leq y\leq z\}$, which is called an {\em order-interval}. Given $y\in A_+$  we let 
\[
\mathrm{face}(y) = \{x\in A_+\colon x\leq \lambda y\mbox{ for some }\lambda\geq 0\}.
\]
Note that $y\sim y'$ if and only if $\mathrm{face}(y) =\mathrm{face}(y')$.  In a Euclidean Jordan algebra $A$ every idempotent $p$ satisfies  
\[
\mathrm{face}(p) \cap [0,u] =[0,p],
\]
by \cite[Lemma 1.39]{AS1}. 

Let us now prove part (a) of Theorem \ref{main}. 
\begin{proof}[Proof of Theorem \ref{main}(a)]
 It follows from Lemma \ref{bijection1}, \ref{cont1}, and \ref{cont2} that $\mathrm{exp}_u\colon \overline{A}^h\to\overline{A_+^\circ}^h$ is a continuous bijection between the compact spaces $\overline{A}^h$ and $\overline{A_+^\circ}^h$. As  $\overline{A_+^\circ}^h$ is Hausdorff, we conclude that $\mathrm{exp}_u$ is a homeomorphism. We know from Theorem \ref{detour} that if $h$ and $h'$ are horofunctions in $\overline{A_+^\circ}^h$, where
\[
h(x) = \max\{\log M(y/x),\log M(z/x^{-1})\} \mbox{\quad  and\quad } 
h'(x) =  \max\{\log M(y'/x),\log M(z'/x^{-1})\}\mbox{\quad for }x\in A_+^\circ,
\] 
then $h$ and $h'$ are in the same part if and only if $y\sim y'$ and $z\sim z'$. Consider the spectral decompositions: $y=\sum_{i\in I} e^{-\alpha_i}p_i $, $y' = \sum_{i\in I'} e^{-\beta_i} q_i$, $z= \sum_{j\in J} e^{-\alpha_j}p_j$, and $z'  = \sum_{j\in j'} e^{-\beta_j} q_j$. If  $y\sim y'$, then $p_I\sim q_{I'}$, where $p_I =\sum_{i\in I} p_i$ and $q_{I'}=\sum_{i\in I'} q_i$. Note that $p_I\sim q_{I'}$  implies $\mathrm{face}(p_I)=\mathrm{face}(q_{I'})$.  As $\mathrm{face}(p_I) \cap [0,u] = [0,p_I]$  by \cite[Lemma 1.39]{AS1}, we get  that $p_I =q_{I'}$.  So $y\sim y'$ implies that $p_I=q_{I'}$. Conversely, if $p_I=q_{I'}$, then $y\sim p_I\sim q_{I'}\sim y'$. Thus, $y\sim y'$ if and only if $p_I=q_{I'}$.  Likewise,  $z\sim z'$ if and only if $p_J=q_{J'}$. Now using \cite[Theorem 4.3]{LP}, we conclude that $\mathrm{exp}_u$ maps parts onto parts. 
  \end{proof}
  
\begin{remark}
It was shown in \cite[Section 4]{LP} that there exists a homeomorphism from the  horofunction compactification of $(A,\|\cdot\|_u)$ onto the closed unit ball in the dual space of $(A,\|\cdot\|_u)$, which maps each part of the horofunction boundary onto a relative open boundary face of the ball. The dual space $(A^*,\|\cdot\|_u^*)$ is a {\em base-norm space}, see \cite[Theorem 1.19]{AS0}. That is to say, it  is an ordered normed vector space with cone $A^*_+=\{\phi\in A^*\colon \phi(x)\geq 0\mbox{ for all }x\in A_+\}$, $A^*_+-A^*_+ = A^*$, and the unit ball of the norm is given by 
\[
B_1^* = \mathrm{conv} (S(A)\cup -S(A)),
\]
where $S(A) =\{\phi\in A_+^*\colon \phi(u)=1\}$ is the {\em state space} of $A$. 

If we identify the  finite dimensional formally real Jordan algebra $A$ with $A^*$ using the inner-product $(x|y) =\mathrm{tr}(x\bullet y)$, we get that  $A^*_+=A_+$, as $A_+^\circ$ is a symmetric cone (see \cite[Proposition III.4.1]{FK}) and $S(A) =\{w\in A_+\colon (u|w)=1\}$.  It was shown in \cite[Theorem 4.4]{ER} that the (closed) boundary faces of the dual ball $B_1^*= \mathrm{conv} (S(A)\cup -S(A))\subset A$ are precisely the sets of the form, 
\begin{equation}\label{ER}
F_{p,q}=\mathrm{conv}\, ( (U_p(A)\cap S(A))\cup(U_q(A)\cap -S(A))),
\end{equation}
where $p$ and $q$ are orthogonal idempotents in $A$. 
\end{remark}

  \section{Variation norm  horofunctions}
  Let $A$ be a finite dimensional formally real Jordan algebra  and $T_u = \{w\in A\colon \mathrm{tr}\,w =0\}$, which is the tangent space of $PA_+^\circ=\{x\in A_+^\circ\colon \det x =1\}$ at the unit $u$.  Consider the variation norm, 
  \[
  |w|_u =M(w/u)-m(w/u) = \mathrm{diam}\, \sigma (w)
  \]
  on $T_u$. In this section we determine the horofunction compactification of the normed space $(T_u,|\cdot|_u)$. 
 
 We start by giving the general form of the horofunctions. 
  \begin{proposition}\label{formofg} 
  If $g\colon T_u\to\mathbb{R}$ is a horofunction of $\overline{T_u}^h$, then there exist a Jordan frame $p_1,\ldots,p_r\in A$, $I,J\subseteq \{1,\ldots,r\}$ disjoint and nonempty, and $\alpha\in \mathbb{R}^{I\cup J}$ with $\min_{i\in I}\alpha_i = 0=\min_{j\in J}\alpha_j$ such that 
  \begin{equation}\label{horog}
  g(v) = \Lambda_{A(p_I)} (-U_{p_I}v-\sum_{i\in I}\alpha_ip_i) + \Lambda_{A(p_J)} (U_{p_J}v-\sum_{j\in J}\alpha_jp_j)
  \end{equation}
  for $v\in T_u$, where $p_I = \sum_{i\in I} p_i$ and $p_J= \sum_{j\in J} p_j$. 
  \end{proposition}
  \begin{proof}
  Suppose that $(w^n)$ in $T_u$ is such that $g_{w^n}\to g\in \partial \overline{T_u}^h$. Then by Lemma \ref{Rieffel} we know that 
  $|w^n|_u = \diam\sigma(w^n) = \Lambda(w^n) +\Lambda(-w^n)\to\infty$. Let    $z^n = w^n-\frac{1}{2}(\Lambda(w^n)-\Lambda(-w^n))u\in A$. Note that for each $v\in A$ we have  that $|v-w^n|_u =  |v-z^n|_u$ by (\ref{seminorm}).   
  Moreover, by construction, $\Lambda(z^n) = \frac{1}{2}(\Lambda(w^n) +\Lambda(-w^n)) = \Lambda(-z^n)$. Let $r_n = \Lambda(z^n)$. 
  
  Using the spectral decomposition we write $z^n = \sum_{i=1}^r \mu_i^n p_i^n$. After taking subsequences we may assume: 
  \begin{enumerate}[(1)]
  \item There exists $I_+\subseteq\{1,\ldots,r\}$ such that for each $n\geq 1$ we have that $\mu_i^n>0$ if and only if $i\in I_+$. 
  \item $p_i^n\to p_i$ for all $i\in\{1,\ldots,r\}$. 
  \end{enumerate}
 For $i\in I_+$ let  $\alpha_i^n =  r_n-\mu_i^n$,  and set $\alpha_i^n = r_n+\mu_I^n$ for $i\not\in I_+$. So, $\alpha_i^n\geq 0$ for all $i$. Taking a further subsequence we may assume that $\alpha_i^n\to\alpha_i\in [0,\infty]$ for all $i$. Let $I=\{i\in I_+\colon \alpha_i<\infty\}$ and $J=\{j\not\in I_+\colon \alpha_j<\infty\}$. 
 
As $\Lambda(z^n)=\Lambda(-z^n)=r_n$, we know that 
 \[
 \min_{i\in I_+} \alpha_i^n  =0 =\min_{j\not\in I_+}\alpha^n_j,
 \]
 and hence $I$ and $J$ are both nonempty. 
 
 It now follows from \cite[Lemma 4.7]{LP} that 
 \[
  \Lambda(-v +z^n -r_nu)  \to  \Lambda_{A(p_I)} (-U_{p_I}v-\sum_{i\in I}\alpha_ip_i)\mbox{\quad and\quad  } 
 \Lambda(v -z^n -r_nu)\to \Lambda_{A(p_J)} (U_{p_J}v-\sum_{j\in J}\alpha_jp_j). 
 \]
 Thus, 
 \begin{eqnarray*}
\lim_{n\to\infty} g_{w^n}(v) & = & \lim_{n\to\infty} |v-w^n|_u - |w^n|_u\\
  & = & \lim_{n\to\infty} |v-z^n|_u - |z^n|_u\\
  & = & \lim_{n\to\infty} \Lambda(-v+z^n) +\Lambda(v-z^n) -2r_n\\
   & = & \lim_{n\to\infty} \Lambda(-v+z^n-r_nu) +\Lambda(v-z^n-r_nu)\\
   & = & \Lambda_{A(p_I)} (-U_{p_I}v-\sum_{i\in I}\alpha_ip_i) +\Lambda_{A(p_J)} (U_{p_J}v-\sum_{j\in J}\alpha_jp_j)
 \end{eqnarray*}
for all $v\in T_u$, which completes the proof.
\end{proof}
The next proposition shows that each function of the form (\ref{horog}) is indeed a horofunction. In fact, we shall see that it is a Busemann point in $\partial\overline{T_u}^h$. 
\begin{proposition}\label{allgs} Suppose that $p_1,\ldots,p_r\in A$ is a Jordan frame, $I,J\subseteq \{1,\ldots,r\}$ are disjoint and nonempty, and $\alpha\in \mathbb{R}^{I\cup J}$ is such that $\min_{i\in I}\alpha_i = 0=\min_{j\in J}\alpha_j$. If 
\begin{equation}\label{zetaomega}
\zeta = \sum_{i\in I}-\alpha_ip_i+\sum_{j\in J}\alpha_jp_j\mbox{\quad and\quad } \omega = p_I -p_J,
\end{equation}
then for  $\xi^t = t\omega+\zeta - \frac{1}{r}\mathrm{tr}(t\omega +\zeta) u\in T_u$ with $t>0$ we have that $g_{\xi^t}\to g$, where $g$ is given by (\ref{horog}), and hence $g$ is a Busemann point. 
\end{proposition}
\begin{proof}
For $t>0$,  $\omega^t = t\omega+\zeta \in A$, and note that $\Lambda(\omega^t) = \Lambda(-\omega^t) = t$ for all $t>0$ large.  
Then by \cite[Lemma 4.7]{LP} we get that 
\[
\lim_{t\to\infty} \Lambda(-v+\omega^t - tu) = \lim_{t\to\infty}\Lambda(-v- \sum_{i\in I}\alpha_ip_i +\sum_{j\in J} (-2t+\alpha_j)p_j+\sum_{k\not\in I\cup J}-tp_k) = \Lambda_{A(p_I)} (-U_{p_I}v-\sum_{i\in I}\alpha_ip_i).
\]
Likewise,
\[
\lim_{t\to\infty} \Lambda(v-\omega^t - tu) = \Lambda_{A(p_J)} (U_{p_J}v-\sum_{j\in J}\alpha_jp_j).
\]
Thus, for $v\in T_u$ we have that 
 \begin{eqnarray*}
\lim_{t\to\infty} g_{\xi^t}(v) & = & \lim_{t\to\infty} |v-\xi^t|_u - |\xi^t|_u\\
  & = & \lim_{t\to\infty} |v-\omega^t|_u - |\omega^t|_u\\
  & = & \lim_{t\to\infty} \Lambda(-v+\xi^t) +\Lambda(v-\xi^t) -2t\\
   & = & \lim_{t\to\infty} \Lambda(-v+\xi^t-tu) +\Lambda(v-\xi^t-tu)\\
   & = & \Lambda_{A(p_I)} (-U_{p_I}v-\sum_{i\in I}\alpha_ip_i) +\Lambda_{A(p_J)} (U_{p_J}v-\sum_{j\in J}\alpha_jp_j), 
 \end{eqnarray*}
 which shows that $g_{\xi^t}\to g$. As $t\mapsto \xi^t$ is a straight-line geodesic, we find that $g$ is Busemann point. 
\end{proof}
By combining Propositions \ref{formofg} and \ref{allgs} we get the following description of the horofunctions in $\overline{T_u}^h$. 
\begin{theorem} \label{horoTu}
The horofunctions of $\overline{T_u}^h$ are precisely the functions $g\colon T_u\to\mathbb{R}$ of the form (\ref{formofg}), and each horofunction is a Busemann point.
\end{theorem}
Let us now analyse the parts and detour distance for $\overline{T_u}^h$. 
\begin{proposition}\label{detourTu}
Let $g,g'\in \partial\overline{T_u}^h$ be two horofunctions, where $g$ is given by (\ref{formofg}) and 
 \begin{equation}\label{horog'}
  g'(v) = \Lambda_{A(q_{I'})} (-U_{q_{I'}}v-\sum_{i\in I'}\beta_iq_i) + \Lambda_{A(q_{J'})} (U_{q_{J'}}v-\sum_{j\in J'}\beta_jq_j)\mbox{\quad for }v\in T_u.
  \end{equation}
  If $p_I = q_{I'}$ and $p_J=q_{J'}$, then $g$ and $g'$ are in the same part and 
  \[
  \delta(g,g') = \Big( \Lambda_{A(p_I)}(a_I-b_{I'}) +\Lambda_{A(p_I)}(b_{I'}-a_I) \Big) 
  + \Big(\Lambda_{A(p_j)}(a_J-b_{J'}) +\Lambda_{A(p_J)}(b_{J'}-a_J)\Big),  
  \]
  where $a_I = \sum_{i\in I}\alpha_ip_i$, $a_J = \sum_{j\in J}\alpha_jp_j$, $b_{I'} = \sum_{i\in I'}\beta_iq_i$ and  $b_{J'} = \sum_{j\in J'}\beta_jq_j$.
 \end{proposition}
\begin{proof}
Let $\zeta$, $\omega$ and $\xi^t$ be as in Proposition \ref{allgs}. Then for all $t>0$ large we have that 
\begin{eqnarray*}
|\xi^t|_u+g'(\xi^t) & = & |\xi^t|_u +\Lambda_{A(q_{I'})} (-U_{q_{I'}}\xi^t-b_{I'}) + \Lambda_{A(q_{J'})} (U_{q_{J'}}\xi^t-b_{J'})\\ 
 & = & 2t + \Lambda_{A(p_I)} (-U_{p_I}\xi^t-b_{I'}) + \Lambda_{A(p_J)} (U_{p_J}\xi^t-b_{J'})\\ 
 & = & \Lambda_{A(p_I)} (tp_I-U_{p_I}\xi^t-b_{I'}) + \Lambda_{A(p_J)} (tp_J+U_{p_J}\xi^t-b_{J'})\\ 
& = & \Lambda_{A(p_I)} (\frac{1}{r}\mathrm{tr}(t\omega +\zeta)p_I+a_I-b_{I'}) + \Lambda_{A(p_J)} (-\frac{1}{r}\mathrm{tr}(t\omega +\zeta)p_J+a_J- b_{J'})\\ 
& = &  \Lambda_{A(p_I)} (a_I -b_{I'}) + \Lambda_{A(p_J)} (a_J-b_{J'}).
\end{eqnarray*}
So from (\ref{detourcost}) and Proposition \ref{allgs}, we conclude that $H(g,g') =  \Lambda_{A(p_I)} (a_I -b_{I'}) + \Lambda_{A(p_J)} (a_J-b_{J'})$. Interchanging the roles of $g$ and $g'$ gives $H(g',g) = \Lambda_{A(p_j)}(a_J-b_{J'}) +\Lambda_{A(p_J)}(b_{J'}-a_J)$, which completes the proof.
\end{proof}
The condition in Proposition \ref{detourTu} characterises the parts in the horofunction boundary as the next proposition shows. 
\begin{proposition}\label{partsTu}
If $g$ and $g'$ are horofunctions in $\overline{T_u}^h$ given by (\ref{horog}) and (\ref{horog'}), respectively, then $g$ and $g'$ are in the same part  if and only if $p_I=q_{I'}$ and $p_J=q_{J'}$. 
\end{proposition}
\begin{proof}
By Proposition \ref{detourTu} it remains to show that $\delta(g,g')=\infty$ if $p_I\neq q_{I'}$ or $p_J\neq q_{J'}$.  Suppose that $p_I\neq q_{I'}$. Then $p_I\nleq q_{I'}$ or $q_{I'}\nleq p_I$. Suppose that $p_I\nleq q_{I'}$. Let $\zeta$, $\omega$ and $\xi^t$ be as in Proposition \ref{allgs}. We will show that $H(g,g')=\infty$ in this case. 

For all $t>0$ large we have that 
\begin{eqnarray*}
|\xi^t|_u+g'(\xi^t) & = & |\xi^t|_u +\Lambda_{A(q_{I'})} (-U_{q_{I'}}\xi^t-b_{I'}) + \Lambda_{A(q_{J'})} (U_{q_{J'}}\xi^t-b_{J'})\\ 
 & = & 2t + \Lambda_{A(q_{I'})} (-U_{q_{I'}}\xi^t-b_{I'}) + \Lambda_{A(q_{J'})} (U_{p_J}\xi^t-b_{J'})\\ 
 & = & \Lambda_{A(q_{I'})} (tq_{I'}-U_{q_{I'}}\xi^t-b_{I'}) + \Lambda_{A(q_{J'})} (tq_{J'}+U_{q_{J'}}\xi^t-b_{J'})\\ 
 & = & \Lambda_{A(q_{I'})}(tq_{I'}-U_{q_{I'}}(t\omega+\zeta)-b_{I'}) + \Lambda_{A(q_{J'})} (tq_{J'}+U_{q_{J'}}(t\omega+\zeta)-b_{J'}).
\end{eqnarray*}
Note that $t\omega+\zeta\leq tp_I$ for all $t>0$ large. So, $U_{q_{I'}}(t\omega+\zeta)\leq tU_{q_{I'}}p_I$ for all $t>0$ large. 
This implies that 
\[
\Lambda_{A(q_{I'})}(tq_{I'}-U_{q_{I'}}(t\omega+\zeta)-b_{I'}) \geq \Lambda_{A(q_{I'})}(t(q_{I'}-U_{q_{I'}}p_I)-b_{I'}) \to\infty
\]
as $t\to\infty$, since $q_{I'}-U_{q_{I'}}p_I>0$ by \cite[Lemma 4.12]{LP}. 

Also note that for all $t>0$ large we have that $t\omega+\zeta\geq -tp_J$, and hence $U_{q_{J'}}(t\omega+\zeta)\geq -tU_{q_{J'}}p_J$. 
As $U_{q_{J'}}p_J \leq U_{q_{J'}}u =q_{J'}$ conclude that $t(q_{J'}-U_{q_{J'}}p_J)\geq 0$ for all $t>0$ large. 
It follows that 
\[
\Lambda_{A(q_{J'})} (tq_{J'}+U_{q_{J'}}(t\omega+\zeta)-b_{J'}) \geq \Lambda_{A(q_{J'})} (t(q_{J'}-U_{q_{J'}}p_J) -b_{J'})\geq 
\Lambda_{A(q_{J'})}(-b_{J'})>-\infty
\]
for all $t>0$ large. Combining this inequality with the previous one and using (\ref{detourcost}), we conclude that $H(g,g')=\infty$. 

For the other cases the result can be shown in the same way. 
\end{proof}

\section{Proof of part (b) of Theorem \ref{main}}
The proof of part (b) of Theorem \ref{main} follows the same steps as the one taken in the proof of part (a). To define the extension we recall the characterisation of the Hilbert distance horofunctions of $PA_+^\circ$ from \cite{LP}. The horofunctions of  $(PA_+^\circ,d_H)$ are all Busemann points and precisely the functions $h\colon PA_+^\circ\to\mathbb{R}$ of the form 
\begin{equation}\label{expg}
h(x) = \log M(y/x)+\log M(z/x^{-1})\mbox{\quad for } x\in PA_+^\circ,
\end{equation}
where $y,z\in\partial A_+$ are such  that $\|y\|_u=\|z\|_u =1$ and $y\bullet z=0$, see \cite[Theorem 5.4]{LP}.
In this case the extension of the exponential map is defined as follows. 
\begin{definition}
 The {\em exponential map}, $\mathrm{exp}_u\colon \overline{T_u}^h\to \overline{PA_+^\circ}^h$, is defined by $\mathrm{exp}_u(v) = e^v$ for $v\in T_u$, and for $g\in\partial \overline{T_u}^h$ given by (\ref{horog}) we let  $\mathrm{exp}_u(g) = h$, where $h$ is given by (\ref{expg}) with 
 $y= \sum_{i\in I}e^{-\alpha_i}p_i$ and $z=\sum_{j\in J}e^{-\alpha_j}p_j$. 
\end{definition} 
Note that the extension is well-defined. Indeed, if $g$ given by (\ref{horog}) is represented as 
\[
g(v) = \Lambda_{A(q_{I'})} (-U_{q_{I'}}v-\sum_{i\in I'}\beta_iq_i) + \Lambda_{A(q_{J'})} (U_{q_{J'}}v-\sum_{j\in J'}\beta_jq_j),
\]
then , as $\delta(g,g)=0$, we get by Propositions \ref{detourTu} and \ref{partsTu}  that $p_I=q_{I'}$, $p_J=q_{J'}$,  
$\sum_{i\in I}\alpha_ip_i = \sum_{i\in I'}\beta_iq_i $ and $ \sum_{j\in J}\alpha_jp_j = \sum_{j\in J'}\beta_jq_j$.  
From Lemma \ref{equal} we deduce that $\sum_{i\in I}e^{-\alpha_i}p_i = \sum_{i\in I'}e^{-\beta_i}q_i $ and $ \sum_{j\in J}e^{-\alpha_j}p_j = \sum_{j\in J'}e^{-\beta_j}q_j$, and hence the extension is well-defined.

We first show that $\mathrm{exp}_u\colon \overline{T_u}^h\to \overline{PA_+^\circ}^h$ is a bijection. 
\begin{lemma}\label{bijH} $\mathrm{exp}_u\colon \overline{T_u}^h\to \overline{PA_+^\circ}^h$ is a bijection which maps $T_u$ onto $PA_+^\circ$, and $\partial \overline{T_u}^h$ onto $\partial \overline{PA_+^\circ}^h$. 
\end{lemma}
\begin{proof}
As $\det \exp_u(x) = e^{\mathrm{tr}\, x} = 1$ for $x\in T_u$, we see that $\mathrm{exp}_u$ is a bijection from $T_u$ onto $PA_+^\circ$.  
It follows from \cite[Theorem 5.4]{LP} and Theorem \ref{horoTu} that $\mathrm{exp}_u$ maps $\partial\overline{T_u}^h$ onto $\partial \overline{PA_+^\circ}^h$. 

To complete the proof it remains to show that if $g,g'\in \partial\overline{T_u}^h$ with $\mathrm{exp}_u(g) =\mathrm{exp}_u(g')$, then $g=g'$. Let $g$ and $g'$ be given by (\ref{horog}) and (\ref{horog'}), respectively. By definition of $\mathrm{exp}_u$ we have that $\mathrm{exp}_u(g) = h$, where 
\[
h(x)= \log M(y/x)+\log M(z/x^{-1})\mbox{\quad for }x\in PA^\circ_+,
\]
with $y= \sum_{i\in I}e^{-\alpha_i}p_i$ and $z=\sum_{j\in J}e^{-\alpha_j}p_j$. 
Likewise, $\mathrm{exp}_u(g') = h'$, where $h'(x)= \log M(y'/x)+\log M(z'/x^{-1})$  for $x\in PA^\circ_+$, with $y'= \sum_{i\in I'}e^{-\beta_i}q_i$ and $z=\sum_{j\in J'}e^{-\beta_j}q_j$.

As $h=h'$, we know by $\delta(h,h') =0$, and hence $y=y'$ and $z=z'$ by \cite[Proposition 5.6]{LP}. By the spectral theorem \cite[Theorems III.1.1 and III.1.2]{FK}, it follows that $\sum_{i\in I}\alpha_ip_i = \sum_{i\in I'}\beta_iq_i$ and $\sum_{j\in J}\alpha_jp_j = \sum_{j\in J'}\beta_jq_j$. 
Moreover,  $p_I =q_{I'}$ and $p_J = q_{J'}$. 
It now follows that $g=g'$, which completes the proof. 
\end{proof}
The proof of  the continuity of the extension of $\mathrm{exp}_u$ is split up into two lemmas. 
\begin{lemma}\label{cont1b}
If $(w^n)$ in $T_u$ converges to $g\in\partial \overline{T_u}^h$, then $(\mathrm{exp}_u(w^n))$ converges to $\mathrm{exp}_u(g)$. 
\end{lemma}
\begin{proof}
Let $(w^n)$ be a sequence in $T_u$ converging to $g\in\partial \overline{T_u}^h$, where $g$ is given by (\ref{horog}). To prove that  
$(\mathrm{exp}_u(w^n))$ converges to $\mathrm{exp}_u(g)=h$, where $h$ is given by (\ref{expg}), we show that each of its subsequences has a subsequence converging to $h$. So let $(\mathrm{exp}_u(w^{n_k}))$ be a subsequence.  As $g$ is a horofunction, it follows from Lemma \ref{Rieffel} that $|w^{n_k}|_u\to \infty$. Set $z^{n_k} = w^{n_k}-\frac{1}{2}(\Lambda(w^{n_k})-\Lambda(-w^{n_k}))u$ and let $r_{n_k}=\Lambda(z^{n_k})$. So, 
$2r_{n_k} = \Lambda(z^{n_k})+ \Lambda(-z^{n_k})=|z^{n_k}|_u=|w^{n_k}|_u$.

Using the spectral decomposition we write $z^{n_k} = \sum_{i=1}^r \lambda^{n_k}_iq^{n_k}_i$. By taking subequences we may assume that there exists $I_+\subseteq \{1,\ldots,r\}$ such that for all $k$ we have that $\lambda_i^{n_k}>0$ if and only if $i\in I_+$. In addition, we may assume that $q^{n_k}_i\to q_i$ for all $i$.  Let $\beta^{n_k}_i = r_{n_k}-\lambda_i^{n_k}\geq 0$ for $i\in I_+$, and $\beta^{n_k}_i = r_{n_k}+\lambda_i^{n_k}\geq 0$ for $i\not\in I_+$. By taking further subsequences we may assume that $\beta^{n_k}_i\to\beta_i\in [0,\infty]$ for all $i$. Let $I'=\{i\in I_+\colon \beta_i<\infty\}$ and $J'=\{j\not\in I_+\colon \beta_j<\infty\}$. 

Note that $\min_{i\in I'}\beta_i =0=\min_{j\in J'}\beta_j$, as $\min_{i\in I_+}\beta^{n_k}_i = 0 = \min_{i\not\in I_+}\beta^{n_k}_i$ for all $k$, so $I'$ and $J'$ are both nonempty. It follows from \cite[Lemma 4.7]{LP} and (\ref{seminorm}) that 
\begin{eqnarray*}
\lim_{k\to\infty} g_{w^{n_k}}(v) & = & \lim_{k\to\infty} |v-w^{n_k}|_u - |w^{n_k}|_u\\
  & = & \lim_{k\to\infty} |v-z^{n_k}|_u - 2r_{n_k}\\
  & = & \lim_{k\to\infty} \Lambda(-v+z^{n_k} -r_{n_k}u)+\Lambda(v-z^{n_k} -r_{n_k}u) \\
  & = & \Lambda_{A(q_{I'})}(-U_{q_{I'}}v-\sum_{i\in I'} \beta_iq_i)+ \Lambda_{A(q_{J'})}(U_{q_{J'}}v-\sum_{j\in J'} \beta_jq_j).\\
\end{eqnarray*}
If we denote the righthand side by $g'(v)$, we find that $g'\colon T_u\to\mathbb{R}$ is a horofunction by Proposition \ref{allgs}. 

Since $g_{w^n}\to g$, we conclude that $g=g'$, and hence $\delta(g,g')=0$. It now follows from Propositions \ref{detourTu} and \ref{partsTu} that 
$p_I =q_{I'}$ and $p_J=q_{J'}$. Moreover, $ \sum_{i\in I}\alpha_ip_i =   \sum_{i\in I'}\beta_iq_i$ and $\sum_{j\in J}\alpha_jp_j = \sum_{j\in J'}\beta_jq_j$. So by Lemma \ref{equal} we get that $y=\sum_{i\in I}e^{-\alpha_i}p_i=   \sum_{i\in I'}e^{-\beta_i}q_i$ and $z=\sum_{j\in J}e^{-\alpha_j}p_j = \sum_{j\in J'}e^{-\beta_j}q_j$. 

It follows that 
\[
e^{-r_{n_k}}\mathrm{exp}_u(z^{n_k}) = \sum_{i=1}^r e^{-(r_{n_k}-\lambda_i^{n_k})}q_i^{n_k}\to \sum_{i\in I'}\beta_iq_i=y
\]
and 
\[
e^{-r_{n_k}}\mathrm{exp}_u(-z^{n_k}) = \sum_{i=1}^r e^{-(r_{n_k}+\lambda_i^{n_k})}q_i^{n_k}\to \sum_{j\in J'}\beta_jq_j=z.
\]

As $\mathrm{exp}_u(a+\lambda u) = e^\lambda \mathrm{exp}_u(a)$ for all $a\in A$ and $\lambda \in\mathbb{R}$, we have that 
\begin{eqnarray*}
\lim_{k\to\infty} h_{\mathrm{exp}_u(w^{n_k})}(x) & = & \lim_{k\to\infty} d_H(x,\mathrm{exp}_u(w^{n_k})) - d_H(u, \mathrm{exp}_u(w^{n_k}))\\
 & = & \lim_{k\to\infty} d_H(x,\mathrm{exp}_u(z^{n_k})) - d_H(u, \mathrm{exp_u}(z^{n_k}))\\
  & = & \lim_{k\to\infty} \log M(\mathrm{exp}_u(z^{n_k})/x)+\log M(x/\mathrm{exp}_u(z^{n_k})) -  \mathrm{diam}\,\sigma(z^{n_k})\\
  & = & \lim_{k\to\infty}  \log M(\mathrm{exp}_u(z^{n_k})/x)+\log M(\mathrm{exp}_u(-z^{n_k})/x^{-1}) -  2r_{n_k}\\
  & = & \lim_{k\to\infty}  \log M(e^{-r_{n_k}}\mathrm{exp}_u(z^{n_k})/x)+\log M(e^{-r_{n_k}}\mathrm{exp}_u(-z^{n_k})/x^{-1})\\
  & = & \log M(y/x)+\log M(z/x^{-1})
  \end{eqnarray*} 
for all $x\in PA_+^\circ$ by continuity of the $M$ functions, see \cite[Lemma 2.2]{LLNW}. This shows that $(\mathrm{exp}(w^{n_k}))$ in $PA_+^\circ$ converges $h$, and hence the proof is complete. 
  \end{proof}
  Next we establish the continuity in the horofunction boundary. 
  \begin{lemma}\label{cont2b}
If $(g_n)$ in $\partial \overline{T_u}^h$ converges to $g\in\partial \overline{T_u}^h$, then $(\mathrm{exp}_u(g_n))$ converges to $\mathrm{exp}_u(g)$. 
\end{lemma}
\begin{proof}
Let $(g_n)$ be a sequence in $\partial \overline{T_u}^h$ converging to $g$, where $g$ is given by (\ref{horog}). So $\mathrm{exp}_u(g)=h$, where is $h$ given by (\ref{expg}). We prove that each subsequence of $(\mathrm{exp}_u(g_n))$ has a convergent subsequence with limit $h$.  Let $(\mathrm{exp}_u(g_{n_k}))$ be a subsequence. 
 By Theorem \ref{formofg} we can write  for $k\geq 1$, 
 \[
 g_{n_k}(v)  = \Lambda_{A(q^k_{I_k})}(-U_{q^k_{I_k}}v -\sum_{i\in I_k}\beta^k_iq^k_i)+\Lambda_{A(q^k_{J_k})}(U_{q^k_{J_k}}v -\sum_{j\in J_k}\beta^k_jq^k_j),
 \]
 where the $q^k_i$ and $q^k_j$'s are orthogonal primitive idempotents, $I_k,J_k\subseteq\{1,\ldots,r\}$ are nonempty and disjoint, and 
 $\min_{i\in I_k}\beta^k_i = 0 = \min_{j\in J_k}\beta_j^k$. 
 
 After taking subsequences we may assume that 
 \begin{enumerate}[(1)]
 \item There exist $I_0,J_0\subseteq \{1,\ldots,r\}$ such that $I_0=I_k$ and $J_0=J_k$ for all $k$.
 \item There exist $i_0\in I_0$ and $j_0\in J_0$ such that $\beta_{i_0}^k=0=\beta_{j_0}^k$ for all $k$.
 \item $\beta_m^k\to\beta_m\in [0,\infty]$ and $q^k_m\to q_m$ for all $m\in I_0\cup J_0$. 
 \end{enumerate}
 Now let $I'=\{i\in I_0\colon \beta_i<\infty\}$ and $J'=\{j\in J_0\colon \beta_j<\infty\}$. So $I'$ and $J'$ are nonempty and  $\min_{i\in I'}\beta_i = 0 = \min_{j\in J'}\beta_j$. 

 As $I'$ is nonempty, it follows from \cite[Lemma 4.7]{LP} that 
\[
\Lambda_{A(q^k_{I_0})}(-U_{q^k_{I_0}}v -\sum_{i\in I_0}\beta^k_iq^k_i)\to \Lambda_{A(q_{I'})}(-U_{q_{I'}}v -\sum_{i\in I'}\beta_iq_i). 
\]
Likewise, as  $J'$ is nonempty, we know that 
\[
\Lambda_{A(q^k_{J_0})}(U_{q^k_{J_0}}v -\sum_{j\in J_0}\beta^k_jq^k_j)
\to
\Lambda_{A(q_{J'})}(U_{q_{J'}}v -\sum_{j\in J'}\beta_jq_j).
\]
Thus, 
 \[
\lim_{k\to\infty} g_{n_k}(v)  = \Lambda_{A(q_{I'})}(-U_{q_{I'}}v -\sum_{i\in I'}\beta_iq_i)+\Lambda_{A(q_{J'})}(U_{q_{J'}}v -\sum_{j\in J'}\beta_jq_j). 
 \]

 So, if we denote the righthand side by $g'(v)$, then $g'\colon A\to\mathbb{R}$ is a horofunction by Proposition \ref{allgs}. 
 As $g_n\to g$, we find that $g=g'$ and hence $\delta(g,g')=0$. 
 
 It now follows from  Propositions \ref{detourTu} and \ref{partsTu} that  $p_I = q_{I'}$ and $p_J=q_{J'}$. Moreover, 
 \[
 \sum_{i\in I} \alpha_ip_i = \sum_{i\in I'} \beta_iq_i \mbox{ \quad and\quad } \sum_{j\in J}\alpha_jp_j =\sum_{j\in J'}\beta_jq_j, 
 \]
 So by Lemma \ref{equal} we get  that 
 $ y= \sum_{i\in I} e^{-\alpha_i}p_i = \sum_{i\in I'} e^{-\beta_i}q_i$ and $z=\sum_{j\in J}e^{-\alpha_j}p_j = \sum_{j\in J'}e^{-\beta_j}q_j$.

 If we now let $\bar{y}^{k} = \sum_{i\in I_0}e^{-\beta^k_i}q^k_i$ and $\bar{z}^k = \sum_{j\in J_0}e^{-\beta^k_j}q^k_j$, then 
 \[
 \lim_{k\to\infty} \bar{y}^{k}  =  \sum_{i\in I'}e^{-\beta_i}q_i =  y\mbox{\quad and\quad } 
 \lim_{k\to\infty} \bar{z}^{k}  =  \sum_{j\in J'}e^{-\beta_j}q_j =  z.
 \]
 Therefore, by continuity of the $M$ function (\cite[Lemma 2.2]{LLNW}), 
 \[
 \lim_{k\to\infty} \mathrm{exp}_u(g_{n_k})(x) = \lim_{k\to\infty} \log M(\bar{y}^k/x)+\log M(\bar{z}^k/x^{-1}) = 
\log M(y/x)+\log M(z/x^{-1}) = h(x),
 \]
 which completes the proof.
\end{proof}
Collecting the results so far we can now easily proof Theorem \ref{main}(b). 
\begin{proof}[Proof of Theorem \ref{main}(b)]
It follows from  Lemmas \ref{bijH}, \ref{cont1b} and \ref{cont2b} that $\mathrm{exp}_u\colon \overline{T_u}^h\to\overline{PA_+^\circ}^h$ is a continuous bijection.  As  $ \overline{T_u}^h$ is compact and  $\overline{PA_+^\circ}^h$ is Hausdorff, we conclude that $\mathrm{exp}_u$ is a homeomorphism. 

Suppose that $g$ and $g'$ are horofunctions  in the same part of $\overline{T_u}^h$, where $g$ is given by (\ref{formofg}) and $g'$ is given by (\ref{horog'}).
It follows from Propositions \ref{detourTu} and \ref{partsTu} that  $p_I = q_{I'}$ and $p_J=q_{J'}$.   By definition $\mathrm{exp}_u(g) = h$, where $h$ is given by (\ref{expg}) with $y= \sum_{i\in I}e^{-\alpha_i}p_i$ and $z=\sum_{j\in J}e^{-\alpha_j}p_j$. Likewise, $\mathrm{exp}_u(g') = h'$, where $h'$ is given by $h'(x) = \log M(y'/x)+\log M(z'/x^{-1})$,  with $y'= \sum_{i\in I'}e^{-\beta_i}q_i$ and $z'=\sum_{j\in J'}e^{-\beta_j}q_j$. 
As  $p_I = q_{I'}$ and  $p_J=q_{J'}$, we have that $y\sim p_I\sim q_{I'}\sim y'$ and $z\sim p_J\sim q_{J'}\sim z'$, and hence $h$ and $h'$ are in the same part by \cite[Proposition 5.6]{LP}. 

Conversely, if $h$ and $h'$ are in the same part, then $y\sim y'$ and $z\sim z'$ by \cite[Proposition 5.6]{LP}, and hence $p_I\sim q_{I'}$ and $p_J\sim q_{J'}$.  This implies that $\mathrm{face}(p_I)=\mathrm{face}(q_{I'})$.  As $\mathrm{face}(p_I) \cap [0,u] = [0,p_I]$  by \cite[Lemma 1.39]{AS1}, we get  that $p_I =q_{I'}$. Likewise we have that $p_J=q_{J'}$. So, if $h$ and $h'$ are in the same part, then $g$ and  $g'$ are in the same part. This completes the proof.
  \end{proof}
  
\begin{remark} 
The horofunction compactification is homeomorphic to the closed unit ball of the dual space $(T_u,|\cdot|_u)^*$. Indeed, it was shown in 
\cite[Section 5]{LP} that there exists a homeomorphism from the  horofunction compactification of $(PA_+^\circ,d_H)$ to the closed unit ball $B_1^*$ of the dual space $(T_u,|\cdot|_u)^*$, which maps parts onto parts. We know from  \cite[Section 5.3]{LP} that the dual space $(T_u,|\cdot|_u)^*$ is given by  $(T_u,\frac{1}{2}\|\cdot\|_u^*)$, where we  use the inner-product $(x|y) =\mathrm{tr}(x\bullet y)$ to identify $T_u^*$ with $T_u$. The  unit ball $B_1^*$ satisfies
\[
B_1^* =2\mathrm{conv}(S(A)\cup -S(A))\cap T_u,
\]
where $S(A) = \{w\in A_+\colon (u|w) = 1\}$ is the state space of $A$. Its (closed) boundary faces are precisely the nonempty sets of the form, 
\[
A_{p,q} = 2\mathrm{conv}\, ( (U_p(A)\cap S(A)\cup(U_q(A)\cap -S(A)))\cap T_u,
\]
where $p$ and $q$ are orthogonal idempotents by \cite[Theorem 4.4]{ER}. 
\end{remark}

\section{Final remarks} 
Symmetric cones $A_+^\circ$ and their projective cones $PA_+^\circ$ are examples of Riemannian symmetric spaces $X=G/K$ of non-compact type. The Finsler metrics of the Thompson distance  and the Hilbert distance are examples of invariant Finsler metrics, which have been chararcterised by Planche \cite{Pl}.  In \cite{HSWW} it  was shown that each generalised Satake compactification of a symmetric space $X=G/K$ of non-compact type can be realised as a horofunction compactification under an invariant Finsler metric, whose restriction to a flat is a norm with polyhedral unit ball.  In \cite[Examples 5.3 and 5.7]{HSWW} the symmetric space $\mathrm{SL}_n(\mathbb{C})/\mathrm{SU}_n$  is consider for $n=3$ and $n=4$. This space corresponds to the projective symmetric cone $P\Pi_n(\mathbb{C})$ consisting of  $n\times n$ positive definite Hermitian matrices with determinant 1.  If we consider the restriction of the unit ball of the Finsler metric $H(I,\cdot) = |\cdot|_I$ to the  flat consisting of diagonal matrices  in $T_I$, then for $n=3$ we get  a hexagon, and for $n=4$ we get a rhombic dodecahedron.  These unit balls correspond  to the invariant Finsler metrics in \cite[Examples 5.3 and 5.7]{HSWW}, where the generalised Satake compactification is considered for the adjoint representation of $\mathrm{SL}_n(\mathbb{C})$. So for $n=3$ and $n=4$, the Hilbert distance on $\mathrm{SL}_n(\mathbb{C})/\mathrm{SU}_n$ realises the generalised Satake compactification with respect to the adjoint representation. It would be interesting to know if this is true for general $n$.

For a symmetric cone $A_+^\circ$ a flat in the tangent space at $u\in A_+^\circ$ is given by $\mathrm{Span}(\{p_1,\ldots,p_r\})$, where $p_1,\ldots,p_r$ is a Jordan frame.  The  restriction to the flat of the Finsler metric $F$  for the Thompson distance  is a polyhedral norm. In fact, if $w\in \mathrm{Span}(\{p_1,\ldots,p_r\})$ has spectral decomposition $w=\sum_{i=1}^r \delta_i p_i$, then  $F(u,w) =\|w\|_u =  \max_{i} |\delta_i|$.  So, the unit ball of $F$ is an $r$-dimensional hypercube. It would be interesting to investigate if the horofunction compactification of $(A_+^\circ,d_T)$ realises a generalised Satake compactification of the symmetric space $A_+^\circ= \mathrm{Aut}(A_+)/\mathrm{Aut}(A_+)_u$, and if so, to which  representation $\tau\colon \mathrm{Aut}(A_+)\to \mathrm{PSL}_n(\mathbb{C})$  it corresponds.  

In view of Theorem \ref{main} it is  also natural to ask if, more generally, for a symmetric space $X=G/K$  of non-compact type with an invariant Finsler metric, the exponential map $\mathrm{exp}_u\colon T_uX\to X$ extends as a homeomorphism between the horofunction compactification of $X$ under the Finsler distance and the horofunction compactification of the tangent space $T_uX$ under the Finsler norm.

\end{document}